\documentclass[12pt]{article}
\usepackage{times}
\usepackage{amsmath}
\usepackage{amssymb}
\usepackage{mathrsfs}
\usepackage{fancyhdr}
\usepackage[a4paper,text={128mm,185mm},centering]{geometry}
\newtheorem{Lemma}{Lemma}[section]
\newtheorem{Remark}[Lemma]{Remark}
\newtheorem{Theorem}[Lemma]{Theorem}

\newtheorem{Corollary}[Lemma]{Corollary}
\newtheorem{Definition}[Lemma]{Definition} 
\newtheorem{Proposition}[Lemma]{Proposition}
\newenvironment{Proof}[1][\unskip]{\par\vskip6pt plus4pt minus2pt%
 \noindent\textbf{Proof #1:\ }}{\hfill\ensuremath{\blacksquare}}
\newenvironment{Reference}[1]{\pfill\textbf{#1} \textit\bgroup}{\egroup\par}
\def\ad{\mathrm{ad}}
\def\Ad{\mathrm{Ad}}
\def\Ass{\mathrm{Ass}}
\def\aut{\mathfrak{aut}}
\def\Aut{\mathrm{Aut}}

\def\Cat{\mathscr{C}}

\def\End{\mathrm{End}}
\def\F{\mathscr{F}}
\def\FF{\mathfrak{F}}
\def\FB{\mathbf{FB}}
\def\g{\mathfrak{g}}
\def\GL{\mathbf{GL}}
\def\GP{\textsc{gp}}
\def\GTS{\mathbf{GTS}}
\def\grp{\mathrm{grp}}
\def\Grp{\mathbf{Grp}}

\def\Hom{\mathrm{Hom}}
\def\hor{\mathrm{hor}}
\def\id{\mathrm{id}}
\def\im{\mathrm{im}}
\def\inf{\mathrm{inf}}
\def\lft{\mathrm{left}}
\def\loc{\mathrm{loc}}
\def\MF{\mathbf{MF}}
\def\Mor{\textsc{Mor}}

\def\Obj{\textsc{Obj}}
\def\P{\mathbb{P}}

\def\PB{\mathbf{PB}}
\def\pfill{\par\vskip6pt plus4pt minus3pt\noindent}
\def\pr{\mathrm{pr}}
\def\PT{\mathbf{PT}}
\def\R{\mathbb{R}}
\def\Rep{\mathbf{Rep}}
\def\S{\mathbb{S}}
\def\Set{\mathbf{Set}}
\def\term{\mathrm{term}}
\def\triv{\mathrm{triv}}

\def\VB{\mathbf{VB}}
\def\Vect{\mathbf{Vect}}
\def\Vert{\mathrm{Vert}}
\def\vtriv{\mathrm{vtriv}}
\title{FUNCTORIALITY OF PRINCIPAL BUNDLES AND CONNECTIONS}
\author{Gustavo Amilcar Salda\~na Moncada \& Gregor Weingart}
\date{April 18th, 2020}
\begin{document}
\maketitle
\parbox{340pt}{\small
 {\bf R\'esum\'e.}
  L'une des plus importantes contributions de la th\'eorie de jauge en
  math\'ematiques est de souligner l'importance des foncteurs d'association.
  En mettant l'accent sur la th\'eorie des cat\'egories nous caract\'erisons
  ces derniers en utilisant deux de leurs propri\'et\'es naturelles. Cette
  caract\'erisation est ensuite utilis\'ee pour \'etablir une \'equivalence
  entre la cat\'egorie des fibr\'es principaux et une certaine cat\'egorie
  de foncteurs. Du point de vue de la g\'eom\'etrie differentielle nous
  d\'ecrivons la particularisation des connexions non--lin\'eaires ou
  d'Ehresmann au cas principal ou lin\'eaire. La propri\'et\'e d'universalit\'e
  des courbures principales, par ailleurs bien connue et largement utilis\'ee,
  est alors employ\'ee afin de caract\'eriser les fibr\'es vectoriels dans
  l'image d'un foncteur d'association donn\'e.
 \\[6pt]
 {\bf Abstract.}
  Perhaps the most important contribution of gauge theory to general
  mathematics is to point out the importance of association functors.
  Emphasizing category theory we characterize association functors by
  two of their natural properties and use this characterization to establish
  an equivalence between the category of principal bundles and a suitably
  defined category of functors. From the point of view of differential
  geometry we detail the specialization of non--linear or Ehresmann to
  principal and linear connections and discuss the widely known and very
  useful universality of principal curvature in order to characterize the
  vector bundles in the image of a given association functor.
 \\[6pt]
 \\[6pt]
 }
\section{Introduction}
\label{intro}
 Principal bundles and their association functors play a fundamental role
 in differential geometry and mathematical physics. Spin structures in
 pseudo--Riemannian geometry are defined right away as special principal
 bundles and the basic tenet of harmonic analysis is that the canonical
 association functor of a pointed homogeneous space is an equivalence of
 categories to the category of homogeneous fiber bundles. Last but not
 least the choice of principal bundle corresponds to the choice of vacuum
 sector in quantum field gauge theories. Nevertheless principal bundles
 tend to obfuscate calculations due to some inevitable arbitrariness, as
 one can see for example in Cartan geometries and in the botched proof of
 Blunder 5.24 in the otherwise excellent reference \cite{spin}. Explicit
 calculations are more easily done using only the existence of association
 functors and the universality of curvature, arguably one of the most
 useful theorems in all of differential geometry.

 En nuce this article brings these reservations against the use of principal
 bundles to a point: We show that a principal bundle $GM$ over a manifold
 $M$ is completely determined by its association functor $\Ass_{GM}$.
 Conversely every functor $\FF$ from a suitable category of model fibers
 to the category of fiber bundles over $M$ satisfying two more or less
 self--evident axioms agrees with the association functor for some
 principal bundle over $M$. Under natural transformations the class of
 all such functors $\FF$ becomes the category $\GTS_M$ of gauge theory
 sectors over $M$, which turns out to be equivalent to the category
 $\PB_M$ of principal bundles over the manifold $M$.

 Category theory is usually not considered to be of particular importance
 to differential geometry, the text books \cite{nodg} and \cite{jlee} as
 well as the article \cite{pt} are notable exceptions to this rule. Besides
 the characterization of associated vector bundles as geometric vector bundles
 in the sense of \cite{sw} the common differential geometer may find little
 of interest in this article. Our main motivation for studying categorical
 properties of principal bundles nevertheless is the need to formulate the
 proper analogue of the concepts of principal bundles and connections in
 non--commutative geometry along the lines of \cite{micho1}, \cite{micho2}
 and \cite{micho3}. Every definition of quantum bundles with quantum
 connections like the one presented in \cite{sal} will necessarily reflect
 functorial properties of principal bundles in classical differential
 geometry.

 \pfill
 In order to provide a more detailed outline of this article we consider
 a Lie group $G$ and the category $\MF_G$ of manifolds $\F$ endowed
 with a smooth left action $\star:\,G\times\F\longrightarrow\F$ under
 smooth $G$--equivariant maps. Every principal $G$--bundle $GM$ over
 a manifold $M$ defines a functor from the category $\MF_G$ of model
 fibers to the category $\FB_M$ of fiber bundles over $M$
 $$
  \Ass_{GM}:\;\;\MF_G\;\longrightarrow\;\FB_M,
  \qquad\F\;\longmapsto\;GM\,\times_G\F\ ,
 $$
 which we may promote to a functor $\Ass^\omega_{GM}:\,\MF_G\longrightarrow
 \FB^\nabla_M$ to the category of fiber bundles with connections in the
 presence of a principal connection $\omega$ on $GM$. This association
 functor maps Cartesian products in $\MF_G$ to Cartesian products in
 $\FB^\nabla_M$ and maps a manifold $\F$ endowed with the tri\-vial
 $G$--action to the trivial fiber bundle $M\times\F$. Our first main
 theorem stipulates that these two properties already characterize
 association functors as the reader can appreciate in Theorem 5.1

 \pfill
 Let us now consider the category $\PB^\nabla_M$ of principal bundles with
 connections over $M$: Objects are triples $(\,G,\,GM,\,\omega\,)$ formed
 by a Lie group $G$ and a principal $G$--bundle $GM$ over $M$ endowed with
 a principal connection $\omega$, while morphisms are tuples $(\,\varphi_\grp,
 \,\varphi\,)$ consisting of a parallel map $\varphi:\,GM\longrightarrow
 \hat GM$ of the underlying principal bundles, which is equivariant over
 the homomorphism $\varphi_\grp:\,G\longrightarrow\hat G$ of Lie groups. The
 canonical factorization of the model homomorphism $\varphi_\grp$ entails
 a factorization of $\varphi$
 $$
  \varphi:\;\;GM
  \;\stackrel{\pr}\longrightarrow\;
  GM/_{\displaystyle\ker^\circ\varphi_\grp}
  \;\stackrel{\overline\pr}\longrightarrow\;
  GM/_{\displaystyle\ker^{\hphantom\circ}\varphi_\grp}
  \;\stackrel{\overline\varphi}\longrightarrow\;
  \hat GM
 $$
 into a parallel projection, a covering and a parallel injective immersion.
 In this sense every morphism in the category $\PB^\nabla_M$ of principal
 bundles with connections over $M$ is a product of just three basic types:
 The removal of a connected isospin subgroup, a covering $\overline\pr$ of
 principal bundles, a generalized spin structure, and a holonomy reduction
 $\overline\varphi$.

 In order to translate this description of generalized spin structures and
 holonomy reductions as basic type morphism between principal bundles into
 a truly functorial description we consider the category $\GTS^\nabla_M$
 of gauge theory sectors with connections over $M$. Its objects are tuples
 $(\,G,\,\FF\,)$ of a Lie group $G$ together with a functor $\FF:\,\MF_G
 \longrightarrow\FB^\nabla_M$ satisfying the assumptions of Theorem
 \ref{natra}. A morphism $(\,\varphi_\grp,\,\Phi\,)$ between two such gauge
 theory sectors is a natural transformation $\Phi:\,\FF\,\circ\,\varphi^*_\grp
 \longrightarrow\hat\FF$ between the functors twisted by the pull back
 $\varphi^*_\grp:\,\MF_{\hat G}\longrightarrow\MF_G$ of the action along
 the homomorphism $\varphi_\grp:\,G\longrightarrow\hat G$ of Lie groups:
 
 \begin{Reference}{%
  Corollary \ref{equc} (Association Functor as Equivalence of Categories)}
 \hfill\break
  For every smooth manifold $M$ the association functor $\Ass$ provides us
  with an equivalence of categories from the category $\PB^\nabla_M$ of
  principal bundles to the category $\GTS^\nabla_M$ of gauge theory sectors
  with connections:
  $$
   \Ass:\;\;\PB^\nabla_M\;\stackrel\simeq\longrightarrow\;\GTS^\nabla_M,
   \qquad(\,G,\,GM,\,\omega\,)\;\longmapsto\;(\,G,\,\Ass^\omega_{GM}\,)\ .
  $$
  In particular two principal $G$--bundles endowed with principal connections
  on $M$ are isomorphic via a parallel, $G$--equivariant homomorphism of
  fiber bundles, if and only if their association functors are naturally
  isomorphic.
 \end{Reference}

 \pfill
 A direct consequence of Corollary \ref{equc} is that association functors
 are not in general full functors, this is they are not surjective on
 morphisms, simply because the action pull back functor $\varphi^*_\grp:\,
 \MF_{\hat G}\longrightarrow\MF_G$ is not a full functor unless the image
 of $G$ in $\hat G$ is dense. In other words there will be more parallel
 smooth homomorphisms of associated fiber bundles than there are
 $G$--equivariant smooth maps between their model fibers unless
 the principal connection $\omega$ has dense holonomy group.

 According to Corollary \ref{equc} a spin structure on an oriented
 pseudo--Rie\-mannian manifold $(\,M,\,g\,)$ can be defined as a functor
 extending the association functor $\MF_{\mathbf{SO}(T)}\longrightarrow
 \FB^\nabla_M$ determined by the oriented orthonormal frame bundle of
 $M$ to a functor $\MF_{\mathbf{Spin}(T)}\longrightarrow\FB^\nabla_M$
 still satisfying the assumptions of Theorem \ref{natra}, the corresponding
 spinor bundle $\$ M$ is simply the image of the irreducible Clifford module
 under the extended functor. A fundamental problem in differential geometry
 related to spin structures is to characterize the vector and fiber
 bundles in the image of a given association functor. A partial answer to
 this problem is given in Proposition 4.6, which opens the way
 to an axiomatic characterization of spinor bundles and highlights the
 universality of principal curvatures.

 \pfill
This paper breaks down into five sections.  Section \ref{fiber} is a
 leisurely introduction to non--linear or Ehresmann connections on fiber
 bundles; we relate their curvature to the commutator of iterated covariant
 derivatives and discuss how non--linear connections specialize to principal
 and linear connections. In Section \ref{pocat} we generalize objects of
 group type in categories with Cartesian products to principal objects.
 Association functors are studied in Section \ref{conaf}, the universality
 of curvature is formulated in Proposition \ref{asscon}. Having proved
 Theorem \ref{natra} in Section \ref{agauge} we define the category of
 gauge theory sectors and establish the equivalence of categories formulated
 in Corollary \ref{equc}. 

 \pfill
 The research project described in this article was inspired by the first
 part of the article \cite{nori} and can be seen as a direct analogue of
 this work in the framework of differential instead of algebraic geometry,
 moreover we address the additional complications brought about by the
 presence of connections.
\section{Fiber Bundles and Non--Linear Connections}
\label{fiber}
 Perhaps the single most important concept in differential geometry is
 the notion of connections or the closely related notion of covariant
 derivatives on a vector or more general on a fiber bundle over a fixed
 manifold $M$. In this section we will modify the standard category $\FB_M$
 of fiber bundles over $M$ to a category more useful for our study, the
 category $\FB^\nabla_M$ of fiber bundles with non--linear connections
 over $M$. Moreover we will discuss principal and linear connections
 in the framework of this category.

 \pfill
 In general a fiber bundle over a manifold $M$ with model fiber manifold
 $\F$ is a manifold $\F M$ endowed with a smooth projection map $\pi:\,
 \F M\longrightarrow M$, which is locally trivializable. The preimage of a point $p\,\in\,M$ under $\pi$ is called
 the fiber of the bundle over $p$, it is a submanifold $\F_pM\,:=\,\pi^{-1}
 (\,p\,)\,\subset\,\F M$ of the total space $\F M$ diffeomorphic to the model
 fiber $\F$. Fiber bundles over $M$ are the objects in the category $\FB_M$,
 morphisms in this category are smooth maps $\varphi:\,\F M\longrightarrow
 \hat\F M$ between the total spaces which commute with the respective
 projections $ \hat\pi\circ \varphi=\pi$ and thus map the fibers of $\F M$ to the fibers of $\hat\F M$. Terminal objects in the category $\FB_M$
 correspond to diffeomorphisms $\pi:\,\hat M\longrightarrow M$ thought of
 as fiber bundles over $M$ with single point fiber. 
 
 \pfill  
 The Cartesian product of
 two fiber bundles $\F M$ and $\hat\F M$ in $\FB_M$ is called the fibered
 product in differential geometry $
  \F M\,\times_M\hat\F M$ and it is defined as the equalizer of $\pi\,\circ\,\pr_L$ and $\hat\pi\,\circ\,\pr_R$ in the
 manifold product $\F M\times\hat\F M$. 
   
\pfill  
 In order to study connections
 in the context of category theory we prefer the following definition:

 \begin{Definition}[Non--linear Connections on Fiber Bundles]
 \hfill\label{cfb}\break
  A non--linear connection on a fiber bundle $\F M$ over a manifold $M$
  is a field $\P^\nabla\,\in\,\Gamma(\,\F M,\,\End\,T\F M\,)$ of projections
  $(\P^\nabla)^2\,=\,\P^\nabla$ on the tangent bundle $T\F M$ such that its
  image distribution equals the vertical foliation:
  $$
   \im\;\Big(\;\P^\nabla_f:\;\;T_f\F M\;\longrightarrow\;T_f\F M\;\Big)
   \;\;\stackrel!=\;\;
   \Vert_f\F M\ .
  $$
 \end{Definition}

 \noindent
 Every non--linear connection $\P^\nabla$ on a fiber bundle $\F M$ allows us
 to define the first order differential
 operator
 \begin{equation}\label{nabla}
  D^\nabla:\;\;\Gamma(\,M,\,TM\,)\;\times\;\Gamma_\loc(\,M,\,\F M\,)
  \;\longrightarrow\;\Gamma_\loc(\,M,\,\Vert\,\F M\,)
 \end{equation}
such that
$$
  (\;D^\nabla_Xf\;)_p
  \;\;:=\;\;
  \Big(\;T_pM\;\stackrel{f_{*,\,p}}\longrightarrow\;T_{f(p)}\F M
  \;\stackrel{\P^\nabla_{f(p)}}\longrightarrow\;\Vert_{f(p)}\F M\;\Big)
  \;X_p\ ,
 $$ 
 which is the non--linear analogue of the classical definition of covariant
 derivatives on vector bundles. Somewhat annoyingly this covariant derivative
 $D^\nabla_Xf$ contains the redundant information $f\,=\,\pi_{\F M}\,\circ\,
 D^\nabla_Xf$, where $\pi_{\F M}$ denotes the vertical tangent bundle
 projection $\Vert\,\F M\longrightarrow\F M$, the simplicity of linear
 and principal connections stems from the fact that we can get rid of
 this redundancy altogether, the reduced covariant derivative $\nabla_Xf$
 captures only the partial derivatives of the section $f$.

 The Nijenhuis or curvature tensor of a non--linear connection $\P^\nabla$
 on a fiber bundle $\F M$ over a manifold $M$ is the horizontal $2$--form
 $R^\nabla$ on the total space $\F M$ of the fiber bundle with values in
 the vertical tangent bundle defined for two arbitrary vector fields
 $X,\,Y$ on $\F M$ by:
 \begin{equation}\label{rcurv}
  R^\nabla(\;X,\;Y\;)
  \;\;=\;\;
  -\;\P^\nabla\,
  [\;(\,\id\,-\,\P^\nabla\,)\,X,\;(\,\id\,-\,\P^\nabla\,)\,Y\;]\ .
 \end{equation}
In particular the curvature $R^\nabla$ measures exactly the
 failure of the horizontal distribution $\ker\,\P^\nabla\,\subseteq\,T\F M$
 associated to $\P^\nabla$ to be integrable. An interpretation of the curvature tensor along classical lines as a
 commutator of covariant derivatives is shown in \cite{SaW}

 \begin{Definition}[Parallel Homomorphisms between Fiber Bundles]
 \hfill\label{ph}\break
  A parallel homomorphism between fiber bundles $\F M$ and $\hat\F M$ over the
  same manifold $M$ endowed with connections $\P^\nabla$ and $\P^{\hat\nabla}$
  respectively is a homomorphism $\varphi:\,\F M\longrightarrow\hat\F M$ of
  fiber bundles  such that the following
  diagram commutes:
  $$
   \begin{picture}(116,55)(0,0)
    \put(  0,46){$T\,\F M$}
    \put( 81,46){$T\,\hat\F M$}
    \put(  0, 0){$T\,\F M$}
    \put( 81, 0){$T\,\hat\F M$}
    \put( 37,50){\vector(+1, 0){42}}
    \put( 55,54){$\scriptstyle\varphi_*$}
    \put( 37, 4){\vector(+1, 0){42}}
    \put( 55, 8){$\scriptstyle\varphi_*$}
    \put( 16,42){\vector( 0,-1){30}}
    \put(  1,25){$\scriptstyle\P^\nabla$}
    \put( 96,42){\vector( 0,-1){30}}
    \put(100,25){$\scriptstyle\P^{\hat\nabla}$}
   \end{picture}
  $$
 \end{Definition}

 \noindent
 The constraint $\hat\pi\,\circ\,\varphi\,=\,\pi$ characterizing homomorphisms
 of fiber bundles in the category $\FB_M$ readily implies $\varphi_*(\,\Vert
 \,\F M\,)\,\subset\,\Vert\,\hat\F M$, hence the homomorphism $\varphi$ of
 fiber bundles is parallel, if and only if $\varphi_*$ maps the horizontal
 distribution of $\F M$ to the horizontal distribution of $\hat\F M$:
 $$
  \varphi\;\textrm{\ parallel}
  \qquad\Longleftrightarrow\qquad
  \varphi_*(\;\ker\;\P^\nabla\;)
  \;\;\subset\;\;
  \ker\;\P^{\hat\nabla}\ .
 $$
 Modifying the category $\FB_M$ we define the category $\FB_M^\nabla$ of
 fiber bundles with connection over $M$, in this category morphisms are
 parallel homomorphisms of fiber bundles.

 In the resulting category terminal objects are still diffeomorphisms
 considered as fiber bundles with single point fibers endowed with the
 zero connection $\P^\nabla\,=\,0$. Besides terminal objects the
 category $\FB^\nabla_M$ has Cartesian products: The fibered product
 $\F M\times_M\hat\F M$ of two fiber bundles $\F M$ and $\hat\F M$
 over $M$ with connections $\P^\nabla$ and $\P^{\hat\nabla}$ carries
 the product connection $(\P^\nabla\oplus\P^{\hat\nabla}):\,
 T(\F M\times_M\hat\F M)\longrightarrow\Vert\,\F M\oplus\Vert\,\hat\F M$
 defined by
 $$
  \left.\frac d{dt}\right|_0(\,f_t,\,\hat f_t\,)\;\longmapsto\;
  \P^\nabla\Big(\,\left.\frac d{dt}\right|_0f_t\,\Big)\,\oplus\,
  \P^{\hat\nabla}\Big(\,\left.\frac d{dt}\right|_0\hat f_t\,\Big)\ ,
 $$
 where $t\longmapsto f_t$ and $t\longmapsto\hat f_t$ are smooth curves in
 $\F M$ and $\hat\F M$ subject to the fibered product constraint $\pi(f_t)
 \,=\,\hat\pi(\hat f_t)$ for all $t$. In light of all these definitions the
 Cartesian product with the base manifold $M$ becomes a functor from the
 category $\MF$ of smooth manifolds to the category $\FB^\nabla_M$
 \begin{equation}\label{pf}
  M\,\times:\;\;
  \MF\;\longrightarrow\;\FB^\nabla_M,\qquad\F\;\longmapsto\;M\,\times\,\F\ ,
 \end{equation}
 because every trivial fiber bundle $M\times\F$ over $M$ comes along with
 the trivial connection $\P^\triv$, namely the projection to the tangent
 bundle of $\F$:
 $$
  T(\,M\times\F\,)\;\stackrel\cong\longrightarrow\;TM\,\times\,T\F
  \;\stackrel{\pi\times\id}\longrightarrow\;M\,\times\,T\F\;\;\cong\;\;
  \Vert(\,M\times\F\,)\ .
 $$
 Evidently the horizontal distribution $TM\times\F\,\subset\,T(M\times\F)$
 is an integrable foliation with leaves $M\,\times\,\{f\}$ for every trivial
 connection $\P^\triv$, in consequence $R^\triv\,=\,0$ vanishes necessarily.
 The product functor $M\times$ defined in equation (\ref{pf}) will feature
 prominently in Sections \ref{pocat} and \ref{agauge}.
 
 \pfill
 Having discussed general non--linear connections on fiber bundles in some
 detail we now  want to specialize to principal and linear connections in
 the second part of this section. Recall first of all that a principal bundle
 modeled on a Lie group $G$ is a smooth fiber bundle $GM$ with model fiber
 $G$ endowed with a smooth right $\rho$, fiber preserving action of $G$ on its
 total space $GM$. Also it is possible to define the affine product $\backslash:
  G M\,\times_M GM \longrightarrow G$.
 
 \pfill 
  The
 automorphism group bundle of a principal bundle $GM$ over a manifold $M$
 is the Lie group bundle $\Aut\;GM$ over $M$ defined by
 \begin{equation}\label{aut}
  \Aut\,GM
  \;\;:=\;\;
  \{\;\;(p,\psi)\;\;|\;\;\psi:\,G_pM\longrightarrow G_pM
  \textrm{\ is $G$--equivariant}\;\;\}
 \end{equation}
 with the bundle projection $\pi_{\Aut\,GM}\,:\,\Aut\,GM\longrightarrow M,
 \,(\,p,\,\psi\,)\longmapsto p$. In mathematical physics the Fr\'echet--Lie
 group $\Gamma(\,M,\,\Aut\,GM\,)$ of all global sections of the automorphism
 bundle is called the gauge group of $GM$.

 The fiber of the Lie group bundle $\Aut\;GM$ over a point $p\,\in\,M$ is a
 Lie group $\Aut_pGM$ isomorphic, although not canonically so, to the original
 group $G$, in particular its Lie algebra $\aut_pGM\,\cong\,\g$ is
 isomorphic to the Lie algebra of $G$. All these Lie algebras assemble
 into a smooth Lie algebra bundle $\aut\,GM$, whose global sections
 $\Gamma(\,M,\,\aut\,GM\,)$ form the Fr\'echet--Lie algebra of the gauge
 group $\Gamma(\,M,\,\Aut\,GM\,)$ of the principal bundle $GM$.

 \begin{Definition}[Principal Connections]
 \hfill\label{pconex}\break
  A principal connection on a principal $G$--bundle $GM$ over a manifold
  $M$ is a non--linear connection $\P^\nabla$ on the fiber bundle $GM$,
  which is invariant under the right action of $G$ on $GM$ in the sense
  that the right translations $R_\gamma:\,GM\longrightarrow GM,\,g
  \longmapsto g\gamma,$ are parallel automorphisms for all $\gamma\,\in\,G$.
 \end{Definition}
 
 \noindent
 In difference to general fiber bundles the vertical tangent bundle
 of a principal bundle $GM$ is trivializable by
$$\vtriv:\Vert\,GM
 \longrightarrow GM\times\g,\;\left.\frac d{dt}\right|_0g_t\longmapsto
 (g_0,\left.\frac d{dt}\right|_0g^{-1}_0g_t).$$ This allows to establish the following well--known result

 \begin{Lemma}[Principal Connection Axiom]
 \hfill\label{cpb}\break
  On every principal $G$--bundle $GM$ the association $\P^\nabla\,
  \longleftrightarrow\,\omega$ characterized by $\omega\,:=\,\vtriv\circ
  \P^\nabla$ induces a bijection between principal connections in the sense
  of Definition \ref{pconex} and $\g$--valued $1$--forms $\omega$ on $GM$
  satisfying the axiom
  $$
   \omega_{g_0\gamma_0}\Big(\;\left.\frac d{dt}\right|_0g_t\,\gamma_t\;\Big)
   \;\;=\;\;
   \Ad_{\gamma^{-1}_0}\;\omega_{g_0}\Big(\;\left.\frac d{dt}\right|_0g_t\;\Big)
   \;+\;\left.\frac d{dt}\right|_0\gamma_0^{-1}\gamma_t
  $$
  for all choices of smooth curves $t\longmapsto g_t$ in $GM$ and curves
  $t\longmapsto\gamma_t$ in $G$.
 \end{Lemma}

 \noindent
 Cartan's Second Structure Equation \cite{gtvp} is a convenient description
 of the image of the composition of the curvature tensor $R^\nabla$ with the
 vertical trivialization $\vtriv$ in terms of the exterior derivative of the
 connection form
 \begin{equation}\label{cstruc}
  \Omega
  \;\;:=\;\;
  \vtriv\,\circ\,R^\nabla
  \;\;\stackrel!=\;\;
  d\omega\;+\;\frac12\,[\,\omega\,\wedge\,\omega\,]\ ,
 \end{equation}
 where $\frac12[\,\omega\wedge\omega\,](X,Y)\,:=\,[\,\omega(X),\,\omega(Y)\,]$.

 \pfill
 The strategy persued for linear connections on vector bundles $VM$ follows
 the model of principal connections closely. The tangent bundle of a vector
 space is canonically trivializable $TV\,\cong\,V\times V$ by taking
 actual derivatives and this becomes via $[\,\Vert\,VM\,]_p\,=\,T(\,V_pM\,)$
 the vertical trivialization $$\vtriv:\Vert\,VM\stackrel\cong\longrightarrow VM
 \oplus VM,\,\left.\frac d{dt}\right|_0v_t\longmapsto v_0\oplus\lim_{t\to0}
 \frac1t(v_t-v_0).$$ This map can be used to project out $\nabla_Xv\,:=\,\vtriv(D^\nabla_Xv)$ the
 redundant information from the covariant derivative $D^\nabla_Xv$ of a
 section $v\,\in\,\Gamma(\,M,\,VM\,)$:
  
 \begin{Definition}[Linear Connections on Vector Bundles]
 \hfill\label{linx}\break
  A linear connection on a vector bundle $VM$ on $M$ is a non--linear
  connection $\P^\nabla$ on $VM$ such that the reduced covariant derivative
  is $\R$--bilinear:
  $$
   \nabla:\;\;\Gamma(\,M,\,TM\,)\;\times\;\Gamma(\,M,\,VM\,)
   \;\longrightarrow\;\Gamma(\,M,\,VM\,)\ .
  $$
 \end{Definition}

 \pfill
 In \cite{SaW} it is showed a proof of the following lemma.
 
 \begin{Lemma}[Characterization of Linear Connections]
 \hfill\label{clcx}\break
  A non--linear connection $\P^\nabla$ on a vector bundle gives rise
  to an $\R$--bilinear covariant derivative $\nabla:\,\Gamma(\,M,\,TM\,)
  \times\Gamma(\,M,\,VM\,)\longrightarrow\Gamma(\,M,\,VM\,)$, if and only
  if the multiplication by every $\lambda\,\in\,\R$ is a parallel endomorphism:
  $$
   \Lambda_\lambda:\;\;VM\;\longrightarrow\;VM,
   \qquad v\;\longmapsto\;\lambda\,v\ .
  $$
 \end{Lemma}
\section{Principal Objects in Categories}
\label{pocat}
 In every category $\Cat$ with terminal objects and Cartesian products
 the notion of a group so fundamental to all of mathematics can be
 generalized to the notion of a group like object in $\Cat$. In this
 section we take this beautiful idea to characterize homogeneous spaces
 with trivial stabilizers, generally known as principal homogeneous
 or affine group spaces, in terms of their structure morphisms. Moreover
 we apply this characterization of affine group spaces to the category
 $\FB^\nabla_M$ of fiber bundles with connections over a manifold $M$
 in order to characterize principal bundles with principal connections.

 \pfill
 A group like object in a category $\Cat$ with terminal objects
 and Cartesian products is an object $G\,\in\,\Obj\,\Cat$ together a choice
 of structure morphisms
 $$
  m:\;\;G\;\times\;G\;\longrightarrow\;G
  \qquad\qquad
  \iota:\;\;G\;\longrightarrow\;G
  \qquad\qquad
  \epsilon:\;\;*\;\longrightarrow\;G
 $$
 in $\Cat$ called the multiplication, the inverse and the neutral element
 respectively with an arbitrary fixed terminal object $*$ such that the
 three diagrams
 $$
  \begin{picture}(120,60)(0,-5)
   \put(  0,41){$G\!\times\!G\!\times\!G$}
   \put( 10, 0){$G\!\times\!G$}
   \put(101, 0){$G$}
   \put( 90,41){$G\!\times\!G$}
   \put( 25,38){\vector(0,-1){27}}
   \put(  1,24){$\scriptstyle m\times\id$}
   \put(105,38){\vector(0,-1){27}}
   \put(108,24){$\scriptstyle m$}
   \put( 42, 4){\vector(1, 0){55}}
   \put( 65, 7){$\scriptstyle m$}
   \put( 51,45){\vector(1, 0){36}}
   \put( 58,48){$\scriptstyle \id\times m$}
  \end{picture}
  \quad
  \begin{picture}(105,60)(0,-5)
   \put(10,41){$G$}
   \put( 0, 0){$G\!\times\!G$}
   \put(85, 0){$G$}
   \put(75,41){$G\!\times\!G$}
   \put(15,38){\vector(0,-1){27}}
   \put(-5,24){$\scriptstyle e\times\id$}
   \put(90,38){\vector(0,-1){27}}
   \put(94,24){$\scriptstyle m$}
   \put(21,40){\vector(2,-1){62}}
   \put(52,27){$\scriptstyle\id$}
   \put(32, 4){\vector(1, 0){49}}
   \put(50, 7){$\scriptstyle m$}
   \put(23,45){\vector(1, 0){47}}
   \put(40,48){$\scriptstyle\id\times e$}
  \end{picture}
  \quad
  \begin{picture}(105,60)(0,-5)
   \put(10,41){$G$}
   \put( 0, 0){$G\!\times\!G$}
   \put(85, 0){$G$}
   \put(75,41){$G\!\times\!G$}
   \put(15,38){\vector(0,-1){27}}
   \put(-4,24){$\scriptstyle\iota\times\id$}
   \put(90,38){\vector(0,-1){27}}
   \put(94,24){$\scriptstyle m$}
   \put(21,40){\vector(2,-1){62}}
   \put(51,27){$\scriptstyle e$}
   \put(32, 4){\vector(1, 0){49}}
   \put(50, 7){$\scriptstyle m$}
   \put(23,45){\vector(1, 0){47}}
   \put(40,48){$\scriptstyle\id\times\iota$}
  \end{picture}
 $$
 all commute, where $e\,=\,\epsilon\circ\term$ equals the composition of
 $\epsilon$ with the terminal morphism $\term:\,G\longrightarrow*$. In the
 category $\Set$ of sets for example the terminal objects are sets with
 exactly one element, hence $\epsilon:\,*\longrightarrow G$ essentially
 corresponds to an element of $G$. In turn the commutative diagrams above
 convert respectively into the associativity, the existence of a neutral
 element and the existence of inverses axiom in the definition of a group.
 In other words group like objects in $\Set$ are just plain groups.

 In categories more complicated than $\Set$ the classification of group like
 objects can be simplified by the use of functors: Every covariant functor
 $\FF:\,\Cat\longrightarrow\hat\Cat$, which maps terminal objects to terminal
 objects and preserves Cartesian products, maps group like objects in the
 category $\Cat$ to group like objects in $\hat\Cat$. The standard forgetful
 functor $\MF\longrightarrow\Set$ from manifolds to sets for examples maps a
 group like object in $\MF$ to a group, albeit a Lie group whose multiplication
 and inverse are smooths maps.

 In the same vein group like objects $G$ in the category $\Grp$ of groups
 carry two different group structures, one for being an object in $\Grp$ and
 the other due to the forgetful functor $\Grp\longrightarrow\Set$. It is a
 rather insightful exercise to verify that these two group structures actually
 agree so that $G$ is necessarily abelian, because its multiplication $m:\,
 G\times G\longrightarrow G$ is a morphism in $\Grp$. In consequence the
 fundamental group $\pi_1(\,G,\,e\,)$ of a topological group $G$ is always
 abelian, because the functor $\pi_1$ maps terminal objects to terminal
 objects and preserves Cartesian products.

 \pfill
 With these examples of the usefulness of functors in combination with a
 categorical definition of groups in mind we want to describe the concept
 of an affine group or principal homogeneous space in terms of category theory.
 Given a group like object $G$ in a category $\Cat$ we define a (right)
 principal $G$--object to be an object $X\,\in\,\Obj\;\Cat$ endowed with
 two structure morphisms
 \begin{equation}\label{stm}
  \rho:\;\;X\;\times\;G\;\longrightarrow\;X
  \qquad\qquad
  \backslash:\;\;X\;\times\;X\;\longrightarrow\;G
 \end{equation}
 in $\Cat$ called action and left division respectively such that the action
 diagrams
 \begin{equation}\label{act}
  \vcenter{\hbox{\begin{picture}(97,60)(0,-5)
   \put(0,41){$X$}
   \put(65,41){$X\!\times\!G$}
   \put(40, 0){$X$}
   \put(12,45){\vector(+1, 0){51}}
   \put(34,48){$\scriptstyle\id\times e$}
   \put(11,38){\vector(+1,-1){27}}
   \put(12,22){$\scriptstyle\id$}
   \put(78,38){\vector(-1,-1){27}}
   \put(66,22){$\scriptstyle\rho$}
  \end{picture}}}
  \qquad\qquad
  \vcenter{\hbox{\begin{picture}(120,60)(0,-5)
   \put(  0,41){$X\!\times\!G\!\times\!G$}
   \put( 10, 0){$X\!\times\!G$}
   \put(101, 0){$X$}
   \put( 90,41){$X\!\times\!G$}
   \put( 26,38){\vector(0,-1){27}}
   \put( -2,22){$\scriptstyle\rho\times\pr_R$}
   \put(105,38){\vector(0,-1){27}}
   \put(109,22){$\scriptstyle\rho$}
   \put( 43, 3){\vector(1, 0){55}}
   \put( 67, 8){$\scriptstyle\rho$}
   \put( 52,45){\vector(1, 0){37}}
   \put( 56,48){$\scriptstyle\pr_L\times m$}
  \end{picture}}}
 \end{equation}
 and the following diagrams encoding simple transitivity all commute:
 \begin{equation}\label{stra}
  \vcenter{\hbox{\begin{picture}(110,60)(0,-5)
   \put(0,41){$X\!\times\!X$}
   \put(80,41){$X\!\times\!G$}
   \put(53, 0){$X$}
   \put(33,45){\vector(+1, 0){47}}
   \put(45,49){$\scriptstyle\pr_L\times\backslash$}
   \put(26,38){\vector(+1,-1){27}}
   \put(25,22){$\scriptstyle\pr_R$}
   \put(88,38){\vector(-1,-1){27}}
   \put(79,22){$\scriptstyle\rho$}
  \end{picture}}}
  \qquad\qquad
  \vcenter{\hbox{\begin{picture}(110,60)(0,-5)
   \put(0,41){$X\!\times\!G$}
   \put(80,41){$X\!\times\!X$}
   \put(53, 0){$G$}
   \put(33,45){\vector(+1, 0){47}}
   \put(45,49){$\scriptstyle\pr_L\times\rho$}
   \put(26,38){\vector(+1,-1){27}}
   \put(25,22){$\scriptstyle\pr_R$}
   \put(88,38){\vector(-1,-1){27}}
   \put(79,22){$\scriptstyle\backslash$}
  \end{picture}}}
 \end{equation}
 In these diagrams $\pr_L$ and $\pr_R$ denote the projections to the leftmost
 and rightmost factor respectively, moreover $m:\,G\times G\longrightarrow G$
 denotes multiplication in $G$ and $e:\,X\longrightarrow G$ the composition
 of $\epsilon$ with the terminal morphism $\term:\,X\longrightarrow*$. Left
 principal objects can be defined in complete analogy simply by switching
 left and right factors.

 Intuitively, a principal object is essentially the group object itself,
 where we have forgotten the neutral element, in fact every group like object
 $G$ in a category $\Cat$ becomes a principal object over itself under the two
 structure morphisms $\rho\,:=\,m$ and $\backslash\,:=\,m\,\circ\,(\iota\times
 \id)$. In the category $\Set$ of sets for example a principal object over
 a group $G$ is a set $X$ endowed with a right action $\rho:\,X\times G
 \longrightarrow X,\,(x,g)\longmapsto xg,$ due to the commutative diagrams
 in (\ref{act}) and an additional application $\backslash:\,X\times X
 \longrightarrow G,\,(x,y)\longmapsto x^{-1}y,$ such that the following
 two axioms are met for all $x,\,y\,\in\,X$ and $g\,\in\,G$
 $$
  x\,(\,x^{-1}\,y\,)
  \;\;=\;\;
  y
  \qquad\qquad
  x^{-1}\,(\,x\,g\,)
  \;\;=\;\;
  g\ ,
 $$
 which reflect the commutative diagrams in (\ref{stra}). In consequence the
 right action of $G$ on $X$ is transitive with trivial stabilizers, once we
 have declared an arbitrary point $x\,\in\,X$ to be the neutral element a
 principal object $X\,\neq\,\emptyset$ becomes indiscernible from the group
 $G$. In linear algebra for example it would be appropriate to define an
 affine space to be a principal object $\mathscr{V}\,\neq\,\emptyset$ under
 the additive group underlying a vector space $V$ over a field $\mathbb{K}$.
  
 \begin{Lemma}[Group Like and Principal Objects in $\FB_M^\nabla$]
 \hfill\label{glo}\break
  For every Lie group $G$ the trivial fiber bundle $M\times G$ over a
  manifold $M$ endowed with the trivial connection and the obvious structure
  morphisms is a group like object in the category $\FB_M^\nabla$ of fiber
  bundles with non--linear connections over $M$. Principal $M\times G$--objects
  are exactly the principal $G$--bundles $GM$ over $M$ endowed with a principal
  connection $\omega$.
 \end{Lemma}

 \begin{Proof}
 The product functor $M\times:\,\MF\longrightarrow\FB^\nabla_M,\,\F
 \longmapsto M\times\F,$ maps of course terminal objects in $\MF$ to
 terminal objects in $\FB^\nabla_M$ and preserves Cartesian products,
 hence it maps the Lie group $G$, a group like object in the category
 $\MF$, to the group like object $M\times G$ in the category $\FB^\nabla_M$.
 Consider now a principal $M\times G$--object in $\FB^\nabla_M$, this is a
 fiber bundle $GM$ over $M$ endowed with a non--linear connection $\nabla$
 and structure homomorphisms:
 $$
  \rho:\;\;GM\;\times_M\,(\,M\times G\,)\;\longrightarrow\;GM
  \qquad
  \backslash:\;\;GM\;\times_M\,GM\;\longrightarrow\;M\times G\ .
 $$
 The obvious diffeomorphism $GM\,\times_M(M\times G)\,\cong\,GM\,\times\,G$
 of fiber bundles provides $GM$ with a fiber preserving right action
 $\rho:\,GM\times G\longrightarrow GM$ such that each fiber $G_pM$ becomes
 a principal $G$--object in the category $\Set$, this is to say that the action
 $\rho$ is simply transitive on fibers. For every $\gamma\,\in\,G$ the element
 morphism $\gamma:\,\{*\}\longrightarrow G$ in the category $\MF$ induces
 moreover a parallel homomorphism in the category $\FB^\nabla_M$ of fiber
 bundles
 $$
  GM
  \;\stackrel{\id\times\term}\longrightarrow\;
  GM\,\times_M(\,M\times\{*\}\,)
  \;\stackrel{\id\times\gamma}\longrightarrow\;
  GM\,\times_M(\,M\times G\,)
  \;\stackrel\rho\longrightarrow\;
  GM\ ,
 $$
 which is just the right multiplication $R_\gamma:\,GM\longrightarrow GM,
 \,g\longmapsto g\gamma$. In turn the non--linear connection $\P^\nabla$
 present on the object $GM$ in $\FB_M^\nabla$ arises from a principal
 connection $\omega$ in the sense of Definition \ref{pconex}.
 \end{Proof}

 \pfill
 A general group like object in the category $\FB_M$ of fiber bundles
 over a mani\-fold $M$ is just a bundle of Lie groups over $M$, a fiber
 bundle $GM$ endowed with the structure of a Lie group on each fiber such
 that the multiplication $m:\,GM\times_MGM\longrightarrow GM$, the inverse
 $\iota:\,GM\longrightarrow GM$ and the neutral element section
 $\epsilon:\,M\longrightarrow GM$ are smooth. Somewhat stronger is the
 concept of a Lie group bundle: A bundle $GM$ of Lie groups, which can
 be trivialized locally by group isomorphisms. Evidently this stronger
 condition is necessary and sufficient for the existence of a non--linear
 connection $\P^\nabla$, under which $GM$ becomes a group like object
 in the category $\FB_M^\nabla$.
\section{Association Functors and Principal Bundles}
\label{conaf}
 Principal bundles are in a sense universal fiber bundles, every given
 principal bundle induces myriad fiber bundles with a large variety of
 model fibers over the same base manifold. The construction of all these
 fiber bundles is functorial in nature and best thought of as a functor,
 the association functor $\Ass^\omega_{GM}$, from a suitably defined category
 $\MF_G$ of model fibers to the category $\FB^\nabla_M$ of fiber bundles with
 connections over a manifold $M$. In this section we study the more important
 properties of association functors, the universality of principal connections
 and their curvature and characterize all vector bundles in the image of a
 fixed association functor.
 
 \pfill
 Besides the categories $\FB_M$ and $\FB^\nabla_M$ of fiber bundles we
 are interested in the category $\MF_G$ of manifolds $\F$ acted upon by
 a fixed Lie group $G$ under smooth $G$--equivariant maps $\varphi:\,
 \F\longrightarrow\hat\F$ as morphisms. Terminal objects are one point
 manifolds $\{*\}$ and Cartesian products in the category $\MF_G$ see
 $G$ acting diagonally on the Cartesian product $\F\times\hat\F$ of the
 manifolds underlying two objects $\F$ and $\hat\F$. Interestingly the
 category $\MF_G$ comes along with a canonical endofunctor, the tangent
 bundle endofunctor
 $$
  T:\;\;\MF_G\;\longrightarrow\;\MF_G,\qquad\F\;\longrightarrow\;T\F\ ,
 $$
 which sends an object $\F\,\in\,\Obj\,\MF_G$ to the tangent bundle of its
 underlying manifold considered as a manifold $T\F$ in its own right, on which
 the Lie group $G$ acts by the differential of its characteristic action
 $\star$ on $\F$:
 $$
  \star_{\scriptscriptstyle T\F}:\;\;G\;\times\;T\F\;\longrightarrow\;T\F,
  \qquad\Big(\,\gamma,\,\left.\frac d{dt}\right|_0f_t\;\Big)\;\longmapsto\;
  \left.\frac d{dt}\right|_0\gamma\,\star\,f_t\ .
 $$
 In order to define the tangent bundle functor on morphisms we observe
 that the differential $\varphi_*:\,T\F\longrightarrow T\hat\F$ of a
 $G$--equivariant map $\varphi:\,\F\longrightarrow\hat\F$ is again
 $G$--equivariant and this observation suggests $T\varphi\,:=\,\varphi_*$.
 It should be noted that the Lie group $G$ provides a distinguished object
 in the category $\MF_G$, namely its Lie algebra $\g\,:=\,T_eG$
 considered just as a manifold endowed with the adjoint representation
 $\Ad:\,G\times\g\longrightarrow\g$. The infinitesimal action links this
 distinguished object to the tangent bundle endofunctor:
 
 \begin{Definition}[Infinitesimal Action]
 \hfill\label{infact}\break
  Consider a smooth left action $\star:\,G\times\F\longrightarrow\F,
  \,(\gamma,f)\longmapsto\gamma\,\star\,f,$ of a Lie group $G$ on a
  smooth manifold $\F$. The infinitesimal action of the Lie algebra $\g$
  of the group $G$ associated to this smooth action $\star$ is defined by
  $$
   \star_\inf\,:\;\;\g\;\times\;\F\;\longrightarrow\;T\F,\qquad
   \Big(\;\left.\frac d{dt}\right|_0\gamma_t,\;f\;\Big)\;\longmapsto\;
   \left.\frac d{dt}\right|_0\gamma_t\;\star\;f\ ,
  $$
  where $t\longmapsto\gamma_t$ with $\gamma_0\,=\,e$ represents the tangent
  vector $\left.\frac d{dt}\right|_0\gamma_t\,\in\,\g$.
 \end{Definition}

 \noindent
 En nuce the infinitesimal action is a natural transformation from
 the endofunctor $\g\times$ to the tangent bundle endofunctor. In fact
 $\star_\inf:\,\g\times\F\longrightarrow T\F$ is $G$--equivariant and
 thus a morphism in $\MF_G$ for all objects $\F$ due to
 $$
  \gamma\,\star_{\scriptscriptstyle T\F}(\,X\,\star_\inf f\,)
  \;\;=\;\;
  \left.\frac d{dt}\right|_0
  (\,\gamma\,\gamma_t\,\gamma^{-1}\,)\,\star\,(\,\gamma\,\star\,f\,)
  \;\;=\;\;
  (\,\Ad_\gamma X\,)\,\star_\inf(\,\gamma\,\star\,f\,)
 $$
 for all $f\,\in\,\F$ and all tangent vectors $X\,=\,\left.\frac d{dt}
 \right|_0\gamma_t$ at $\gamma_0\,=\,e$, moreover $\star_\inf$ intertwines
 with the differential $\varphi_*$ of every $G$--equivariant smooth map
 $\varphi:\,\F\longrightarrow\hat\F$ in the identity $\varphi_*(X\star_\inf f)
 \,=\,X\star_\inf\varphi(f)$.

 \begin{Definition}[Association Functor]
 \hfill\label{defass}\break
  Consider a Lie group $G$ and a principal $G$--bundle $GM$ over a
  manifold $M$. Every smooth action $\star:\,G\times\F\longrightarrow\F$
  of the group $G$ on a manifold $\F$ extends to a free and smooth right
  action of the group $G$ on the Cartesian product $GM\times\F$ via
  $(\,g,\,f\,)\,\star\,\gamma\,:=\,(\,g\gamma,\,\gamma^{-1}\star f\,)$.
  The quotient of $GM\times\F$ by this free action is a fiber bundle
  over $M$ with model fiber $\F$
  $$
   \Ass_{GM}(\,\F\,)
   \;\;=\;\;
   GM\,\times_G\F
   \;\;:=\;\;
   (\,GM\,\times\,\F\,)/_{\displaystyle G}
  $$
  called the fiber bundle associated to $GM$ and $\F\,\in\,\Obj\;\MF$.
  Every $G$--equivariant map $\varphi:\,\F\longrightarrow\hat\F$
  induces a homomorphism of fiber bundles
  $$
   \Ass_{GM}(\,\varphi\,):\;\;
   GM\,\times_G\F\;\longrightarrow\;GM\times_G\hat\F,
   \qquad[\,g,f\,]\;\longmapsto\;[\,g,\varphi(f)\,]\ ,
  $$
  which is well--defined in terms of representatives $(g,f)$ of
  the equivalence class $[g,f]$. In other words $\Ass_{GM}:\,\MF_G
  \longrightarrow\FB_M,\,\F\longmapsto GM\times_G\F,$ is a functor
  from $\MF_G$ to the category $\FB_M$ of fiber bundles over $M$.
 \end{Definition}

 \pfill
 Recall now that the each of the categories $\MF_G$ and $\FB_M$ has a
 canonical endofunctor associated with it, namely the tangent bundle
 endofunctor $T$ for the category $\MF_G$ of manifolds with $G$--action
 and the vertical tangent bundle functor $\Vert$ for the category $\FB_M$.
 Considered as a fiber bundle over $M$ the vertical tangent bundle has
 fiber $[\,\Vert\,\F M\,]_p\,=\,T[\,\F_pM\,]$ over every $p\,\in\,M$
 and so we may suspect that the following diagram commutes
 \begin{equation}\label{tcd}
  \vcenter{\hbox{\begin{picture}(130,64)(0,0)
   \put(  0,50){$\MF_G$}
   \put(  0, 0){$\MF_G$}
   \put( 95,50){$\FB_M$}
   \put( 95, 0){$\FB_M$\ .}
   \put( 32,53){\vector(+1, 0){59}}
   \put( 45,57){$\Ass_{GM}$}
   \put( 32, 3){\vector(+1, 0){59}}
   \put( 45, 7){$\Ass_{GM}$}
   \put( 13,46){\vector( 0,-1){35}}
   \put(  0,25){$T$}
   \put(106,46){\vector( 0,-1){35}}
   \put(111,25){$\Vert$}
  \end{picture}}}
 \end{equation}
 up to a natural isomorphism $\Vert(GM\times_G\F)\stackrel\cong
 \longrightarrow GM\times_GT\F$ given by:
 \begin{equation}\label{vtass}
  \left.\frac d{dt}\right|_0\Big[\;g_t,\;f_t\;\Big]
  \;\longmapsto\;\Big[\;g_0,\;\left.\frac d{dt}\right|_0(\,g_0^{-1}\,g_t\,)
  \;\star\;f_t\;\Big]\ .
 \end{equation}
 Of course this isomorphism is motivated by $[\,g_t,f_t\,]\,=\,
 [\,g_0,(g_0^{-1}g_t)\star f_t\,]$, whenever the representative curve
 $t\longmapsto[\,g_t,\,f_t\,]$ for a vertical tangent vector to
 $GM\times_G\F$ has been chosen such that $g_t$ stays in the
 fiber of $g_0$ for all $t$.

 \begin{Remark}[Action of Automorphism Group Bundle]
 \hfill\label{agg}\break
  The automorphism group bundle of a principal bundle $GM$ acts naturally
  $\star:\,\Aut\,GM\times_M\F M\longrightarrow\F M$ on every fiber bundle
  $\F M\,:=\,GM\times_G\F$ associated to $GM$ and an object $\F\,\in\,
  \Obj\;\MF_G$ by means of
  $$
   (\,p,\,\psi\,)\;\star\;[\,g,\,f\,]
   \;\;:=\;\;
   [\,\psi(\,g\,),\,f\,]
  $$
  for all $(\,p,\,\psi\,)\,\in\,\Aut_pGM$ and $[\,g,\,f\,]\,\in\,
  \F_{\pi(g)}M$ in the fibers of $\Aut\,GM$ and $\F M$ over the same
  point $p\,=\,\pi(g)$ of the base manifold $M$.
 \end{Remark}

 \pfill
 In concrete examples the automorphism group bundle $\Aut\,GM$ is usually
 more readily identified than the principal bundle $GM$ itself due to its
 omnipresent action on associated fiber bundles. Consider the orthonormal
 frame bundle of a pseudo--Riemannian manifold $(\,M,\,g\,)$ for example
 $$
  \mathbf{O}(\,M,\,g\,)
  \;\;:=\;\;
  \{\;\;(\,p,\,F\,)\;\;|\;\;p\,\in\,M\textrm{\ and\ }
  F:\;T\longrightarrow T_pM\textrm{\ isometry}\;\;\}\ ,
 $$
 where $T$ is a pseudo--euclidean model vector space of the correct signature
 and $\mathbf{O}(\,T\,)$ acts from the right by precomposition $(\,p,\,F\,)
 \,\gamma\,=\,(\,p,\,F\circ\gamma\,)$. The automorphism group bundle of the
 orthonormal frame bundle $\mathbf{O}(\,M,\,g\,)$ equals the Lie group bundle
 of all infinitesimal isometries of tangent spaces
 $$
  \mathbf{O}(\,TM,\,g\,)
  \;\;:=\;\;
  \{\;\;(\,p,\,\psi\,)\;\;|\;\;\psi:\,T_pM\longrightarrow T_pM
  \textrm{\ isometry}\;\;\}
 $$
 acting by postcomposition $(\,p,\,\psi\,)\star(\,p,\,F\,)\,=\,(\,p,\,
 \psi\circ F\,)$; it just as well acts on the tangent bundle $TM$ and all
 kinds of the tensor bundles etc.

 For a general principal bundle $GM$ we can use the same idea to identify
 the automorphism group bundle $\Aut\,GM$ as a Lie group bundle over $M$
 with the image of a group object in the category $\MF_G$. Letting $G$ act on
 itself by conjugation $\star:\,G\times G\longrightarrow G,\,(\gamma,g)
 \longmapsto\gamma g\gamma^{-1},$ we obtain in fact a group object
 $G^\ad\,\in\,\Obj\,\MF_G$, whose image under the association functor
 is a Lie group bundle $\Ass_{GM}(\,G^\ad\,)$ over $M$ acting
 $G$--equivariantly on $GM$ by
 \begin{equation}\label{giso}
  \Ass_{GM}(\,G^\ad\,)\,\times_MGM\;\longrightarrow\;GM,
  \quad\;(\,[\,g,\gamma\,],\hat g\,)\;\longmapsto\;g\gamma
  (\,g^{-1}\hat g\,)
 \end{equation}
 for all $\gamma\,\in\,G$ and all $g,\,\hat g\,\in\,GM$ in the same fiber.
 In particular $\Aut_pGM$ is isomorphic, but not naturally so, to the Lie
 group $G$ in every $p\,\in\,M$.

 Under this identification $\Ass_{GM}(\,G^\ad\,)\,=\,\Aut\,GM$ of Lie group
 bundles the natural action of $\Aut\,GM$ on associated fiber bundles
 $GM\times_G\F$ pointed out in Remark \ref{agg} becomes the functorial
 extension of the original action $\star$ considered as a $G$--equivariant
 smooth map $\star:\,G^\ad\times\F\longrightarrow\F$. In the same vein the
 functor $\Ass_{GM}$ converts the infinitesimal action of Definition
 \ref{infact} considered as a $G$--equivariant map $\star_\inf:\,\g\times
 \F\longrightarrow T\F$ into
 $$
  \Ass_{GM}(\,\star_\inf\,):\;\;
  (\,GM\,\times_G\g\,)\;\times_M\,(\,GM\,\times_G\F\,)
  \;\longrightarrow\;(\,GM\,\times_GT\F\,)\ ,
 $$
 which in turn becomes the infinitesimal action associated to Remark
 \ref{agg}:
 $$
  \star_\inf:\;\;\aut\;GM\,\times_M\,(\,GM\,\times_G\F\,)
  \;\longrightarrow\;\Vert(\,GM\,\times_G\F\,)\ .
 $$
 Before we proceed to prove the universality of principal connections and
 their curvature we want to digress a little to discuss the gauge principle,
 a fundamental principle in the study of principal bundles allowing us to
 translate calculations on $GM$ to statements about $M$. In its most basic
 formulation the gauge principle is the assertion that we have a canonical
 bijection
 \begin{equation}\label{gp}
  [\,\Omega^\bullet_\hor(\,GM,\,V\,)\,]^G\;\stackrel\cong\longrightarrow\;
  \Omega^\bullet(\,M,\,GM\times_GV\,),\qquad\eta\;\longmapsto\;\GP[\,\eta\,]
 \end{equation}
 between the horizontal differential forms $\eta\,\in\,
 \Omega^\bullet_\hor(\,GM,\,V\,)$ on $GM$ with values in some representation
 $V$ of $G$ satisfying $R_\gamma^*\eta\,=\,\gamma\star\eta$ for all
 $\gamma\,\in\,G$ and general differential forms on the base manifold
 $M$ with values in the associated vector bundle $GM\times_GV$.
 Explicitly this gauge principle reads
 $$
  \GP[\,\eta\,]_p(\,X_1,\,\ldots,\,X_r\,)
  \;\;:=\;\;
  [\,g,\,\eta_g(\,\tilde X_1,\,\ldots,\,\tilde X_r\,)\,]
 $$
 for arbitrary lifts $\tilde X_1,\,\ldots,\,\tilde X_r\,\in\,T_gGM$ of the
 argument tangent vectors $X_1,\,\ldots,\,X_r\,\in\,T_pM$ to an arbitrary
 point $g\,\in\,G_pM$ in the fiber over $p\,\in\,M$. Due to horizontality the
 resulting differential form $\GP[\,\eta\,]$ does not depend on the choice
 of lifts and the assumption $R_\gamma^*\eta\,=\,\gamma\star\eta$ ensures
 that $\GP[\,\eta\,]$ does not depend on the choice of $g\,\in\,G_pM$ either.
 The gauge principle converts the curvature $2$--form $\Omega\,\in\,
 \Omega^2_\hor(\,GM,\,\g\,)$ of Cartan's Second Structure Equation
 (\ref{cstruc}) into a $2$--form on $M$ with values in $\aut\,GM$:
 \begin{equation}\label{c2form}
  R^\omega
  \;\;:=\;\;
  \GP[\;\Omega\;]
  \;\;\in\;\;
  \Omega^2(\;M,\,\aut\,GM\;)
 \end{equation}
 
 \begin{Proposition}[Universality of Principal Curvature]
 \hfill\label{asscon}\break
  Every choice of a principal connection $\omega$ on a principal $G$--bundle
  $GM$ allows us to promote the association functor $\Ass_{GM}:\,\MF_G
  \longrightarrow\FB_M$ to a functor to the category of fiber bundles
  over $M$ with non--linear connections:
  $$
   \Ass^\omega_{GM}:\;\;\MF_G\;\longrightarrow\;\FB^\nabla_M,
   \qquad\F\;\longmapsto\;GM\,\times_G\F\ .
  $$
  In other words $\omega$ induces a natural connection $\nabla$ on
  $GM\times_G\F$ for every $G$--manifold $\F\,\in\,\Obj\,\MF_G$. The
  curvature $R^\nabla$ of this induced connection is determined by
  the infinitesimal action of the Lie algebra bundle $\aut\,GM$
  $$
   \star_\inf:\;\;\aut\;GM\,\times_M\,(\,GM\times_G\F\,)
   \;\longrightarrow\;\Vert(\,GM\times_G\F\,)
  $$
  and the $2$--form $R^\omega\,\in\,\Omega^2(\,M,\,\aut\,GM\,)$. More
  precisely for all local sections $f\,\in\,\Gamma_\loc(\,M,GM\times_G\F\,)$
  and all $X,\,Y\,\in\,\Gamma(\,M,TM\,)$ we find:
  $$
   R^\nabla_{X,\,Y}f
   \;\;=\;\;
   R^\omega(\,X,\,Y\,)\,\star_\inf\,f\ .
  $$
 \end{Proposition}

 \begin{Proof}
 By definition $GM\times_G\F$ is the quotient of the Cartesian product
 $GM\times\F$ by a free right action of the Lie group $G$. In turn the
 canonical projection $\pr:\,GM\times\F\longrightarrow GM\times_G\F$
 defines a tower of fiber bundles
 \begin{equation}\label{tower}
  \vcenter{\hbox{\begin{picture}(55,90)
   \put( 2,80){$GM\,\times\,\F$}
   \put( 0,40){$GM\,\times_G\F$}
   \put(23, 0){$M$}
   \put(30,76){\vector( 0,-1){24}}
   \put(35,62){$\pr$}
   \put(30,36){\vector( 0,-1){24}}
   \put(35,22){$\pi$}
  \end{picture}}}
 \end{equation}
 over $M$, which becomes $U\times(G\times\F)\stackrel\pr\longrightarrow
 U\times\F\stackrel\pi\longrightarrow U$ in a local equivariant trivialization
 of $GM$. The central idea of the proof is to
 choose the connection $\P^\nabla$ on $GM\times_G\F$ such that $\pr$ is
 parallel with respect to the product $\P^\omega\times\P^\triv$ of the
 principal connection $\omega$ on $GM$ and the trivial connection
 $\P^\triv$ on $M\times\F$.

 For this purpose we consider a curve $t\longmapsto(\,g_t,f_t\,)$ in
 $GM\times\F$ and choose a curve $t\longmapsto\gamma_t$ in $G$ with
 $\gamma_0\,=\,e$ representing the tangent vector $\left.\frac d{dt}
 \right|_0\gamma_t\,=\,\omega_{g_0}(\,\left.\frac d{dt}\right|_0g_t\,)
 \,\in\,\g$. The Principal Connection Axiom \ref{cpb} ensures
 $$
  \omega_{g_0e}\Big(\;\left.\frac d{dt}\right|_0g_t\,\gamma_t^{-1}\;\Big)
  \;\;=\;\;
  \Ad^{-1}_e\;\omega_{g_0}\Big(\;\left.\frac d{dt}\right|_0g_t\;\Big)
  \;+\;\left.\frac d{dt}\right|_0e^{-1}\,\gamma^{-1}_t
  \;\;=\;\;
  0
 $$
 and so $t\longmapsto g_t\gamma^{-1}_t$ represents a horizontal tangent
 vector. In turn
 \begin{eqnarray*}
  \lefteqn{(\,\P^\omega\times\P^\triv\,)
  \Big(\;\left.\frac d{dt}\right|_0(\,g_t,\,f_t\,)\;\Big)}
  \quad
  &&
  \\[-5pt]
  &=&
  (\,\P^\omega\times\P^\triv\,)\Big(\;\left.\frac d{dt}\right|_0
  (\,g_t\gamma_t^{-1}\gamma_0,\,f_0\,)\;+\;\left.\frac d{dt}\right|_0
  (\,g_0\,\gamma_0^{-1}\,\gamma_t,\,f_t\,)\;\Big)
  \\[-5pt]
  &=&
  \left.\frac d{dt}\right|_0(\,g_0\,\gamma_t,\,f_t\,)\ ,
 \end{eqnarray*}
 because the first summand is horizontal and the second vertical in
 $GM\times\F$. Projecting this identity to equivalence classes in
 $GM\times_G\F$ we find
 \begin{eqnarray*}
  \P^\nabla\Big(\;\left.\frac d{dt}\right|_0[\,g_t,\,f_t\,]\;\Big)
  &:=&
  \left.\frac d{dt}\right|_0[\,g_0,\,\gamma_t\,\star\,f_t\,]
  \\
  &=&
  \left[\;g_0,\;\left.\frac d{dt}\right|_0f_t\;+\;\omega_{g_0}\Big(\,
  \left.\frac d{dt}\right|_0g_t\,\Big)\,\star_\inf\,f_0\;\right]
 \end{eqnarray*}
 due to the Definition \ref{infact} of the infinitesimal action and the
 choice of the curve $t\longmapsto\gamma_t$. In light of the isomorphism
 (\ref{vtass}) the right hand side denotes a vertical tangent vector to
 $GM\times_G\F$ and so the latter formula defines a non--linear connection
 $\P^\nabla$ on the fiber bundle $GM\times_G\F$.

 With respect to this non--linear connection $\P^\nabla$ the canonical
 projection $\pr:\,GM\times\F\longrightarrow GM\times_G\F$ is parallel,
 because it maps horizontal tangent vectors $\left.\frac d{dt}\right|_0
 [\,g_t,f_0\,]$ with $\omega_{g_0}(\left.\frac d{dt}\right|_0g_t)\,=\,0$
 to horizontal vectors. The construction of $\P^\nabla$ is natural in the
 category $\MF_G$ as well: The functorial extension $\Ass_{GM}(\,\varphi\,):
 \,GM\times_G\F\longrightarrow GM\times_G\hat\F,\,[\,g,f\,]\longmapsto[\,g,
 \varphi(f)\,],$ of every $G$--equivariant smooth map $\varphi:\,
 \F\longrightarrow\hat\F$ is parallel
 \begin{eqnarray*}
  \P^{\hat\nabla}
  \Big(\,\left.\frac d{dt}\right|_0[\,g_t,\varphi(f_t)\,]\,\Big)
  &=&
  \left[\,g_0,\,\left.\frac d{dt}\right|_0\varphi(f_t)\,+\,
  \omega_{g_0}\Big(\left.\frac d{dt}\right|_0g_t\Big)\star_\inf
  \varphi(f_0)\,\right]
  \\
  &=&
  \Ass_{GM}(\,\varphi_*\,)\;\P^{\nabla}
  \Big(\,\left.\frac d{dt}\right|_0[\,g_t,\,f_t\,]\,\Big)
 \end{eqnarray*}
 due to the infinitesimal equivariance $X\star_\inf\varphi(f)\,=\,
 \varphi_*(X\star_\inf f)$. In order to calculate the curvature of the
 connection $\P^\nabla$ we use the fact that in a tower of fiber bundles
 like (\ref{tower}) with a parallel submersion $\pr$ the curvature of
 the image connection $\P^\nabla$ is just the image of the preimage
 connection $\P^\omega\times\P^\triv$ under the differential $\pr_*$.
 Using arbitrary lifts $\tilde X,\,\tilde Y\,\in\,T_gGM$ of tangent
 vectors $X,\,Y\,\in\,T_pM$ to a point $g\,\in\,G_pM$ we calculate
 in this way
 \begin{eqnarray*}
  R^\nabla_{[\,g,\,f]}(\,\tilde X,\,\tilde Y\,)
  &=&
  \pr_{*,\,(\,g,\,f\,)}R^{\omega\times\triv}_{(\,g,\,f\,)}
  (\,\tilde X,\,\tilde Y\,)
  \\
  &=&
  \left.\frac d{dt}\right|_0\;\left[\;g\,
  \exp\Big(\,t\,\Omega_g(\,\tilde X,\,\tilde Y\,)\,\Big),\;f\;\right]
  \\
  &=&
  [\,g,\,\Omega_g(\tilde X,\tilde Y)\,]
  \,\star_\inf\,[\,g,f\,]
  \;\;=\;\;
  R^\omega_p(X,Y)\,\star_\inf\,[\,g,f\,]\ ,
 \end{eqnarray*}
 where $R^\omega\,:=\,\GP[\,\Omega\,]\,\in\,\Omega^2(\,M,\aut\,GM\,)$
 is the $2$--form with values in $\aut\,GM$ the gauge principle (\ref{gp})
 associates to $\Omega\,:=\,d\omega\,+\,\frac12[\,\omega\wedge\omega\,]$.
 Formulated in terms of local sections $f\,\in\,\Gamma_\loc(\,M,\,
 GM\times_G\F\,)$ the latter identity becomes $R^\nabla_{X,\,Y}f\,=\,
 R^\omega(X,Y)\,\star_\inf f$.
 \end{Proof}
 
 \pfill
 One of the most important properties of association functors is that they
 intertwine the actions of smooth functors on the categories $\Rep_G$ and
 $\VB_M$. A smooth functor is an endofunctor $\S:\,\Vect^\times_\R
 \longrightarrow\Vect^\times_\R$ of the category of finite dimensional
 vector spaces under linear isomorphisms such that
 $$
  \Mor_{\Vect^\times_\R}(\,V,\,V\,)\;\longrightarrow\;\Mor_{\Vect^\times_\R}
  (\,\S\,V,\,\S\,V\,),\qquad \varphi\;\longmapsto\;\S(\,\varphi\,),
 $$
 is a smooth map between the smooth manifolds $\Mor_{\Vect^\times_\R}(V,V)
 \,=\,\GL\,V$ and $\Mor_{\Vect^\times_\R}(\S\,V,\S\,V)$ for every finite
 dimensional vector space $V$ over $\R$. Smooth functors extend naturally
 to endofunctors of the category $\Rep_G$ of representations $V$ of a Lie
 group $G$ by letting $G$ act on $\S\,V$ via:
 $$
  \star\;:\;\;G\;\times\;\S\,V\;\longrightarrow\;\S\,V,\qquad
  (\,\gamma,\,s\,)\;\longmapsto\;\S(\;\gamma\,\star:\;
  V\,\stackrel\cong\longrightarrow\,V\;)\;s\ .
 $$
 This extension to representations makes the classification of smooth
 functors an exercise in the representation theory of general linear
 groups: Every smooth functor is naturally isomorphic $\S\,\cong\,
 \S_1\oplus\ldots\oplus\S_r$ to a finite direct sum of Schur functors
 $\S_1,\,\ldots,\,\S_r$ twisted by density lines \cite{fh}.

 In the same vein every smooth functor $\S$ extends naturally to an
 endo\-functor of the category $\VB^\nabla_M$ of vector bundles with
 connections over a manifold $M$. The smoothness of $\S$ allows us
 to define a differentiable structure on the disjoint union of vector
 spaces obtained by applying $\S$ fiberwise
 $$
  \S\;VM
  \;\;:=\;\;
  \bigcup_{p\,\in\,M}\S\;V_pM
 $$
 to obtain a new vector bundle $\S\,VM$ over $M$; every connection
 $\nabla$ on the original vector bundle $VM$ extends naturally to a
 connection $\nabla^\S$ on $\S\,VM$ by the requirement that parallel
 transport with respect to this connection along an arbitrary curve
 $t\longmapsto p_t$ in the manifold $M$ is simply the image
 $$
  \Big(\;\PT^{\nabla^\S}_t:\;
  \S\;V_{p_0}M\,\stackrel\cong\longrightarrow\,\S\;V_{p_t}M\;\Big)
  \;\;=\;\;
  \S\;\Big(\;\PT^\nabla_t:\;
  V_{p_0}M\,\stackrel\cong\longrightarrow\,V_{p_t}M\;\Big)
 $$
 of parallel transport with respect to $\nabla$ under the functor
 $\S$. Because parallel transport in associated vector bundles is
 essentially the image of parallel transport in the principal bundle
 $GM$ itself, every association functor $\Ass^\omega_{GM}$ intertwines
 the two extensions of a smooth functor $\S$ to the categories $\Rep_G$
 of representations and $\VB^\nabla_M$ of vector bundles with connections:
 \begin{equation}\label{sq}
  \vcenter{\hbox{\begin{picture}(130,62)(0,0)
   \put(  0,50){$\Rep_G$}
   \put(102,50){$\VB^\nabla_M$}
   \put(  0, 0){$\Rep_G$}
   \put(102, 0){$\VB^\nabla_M$\ .}
   \put( 33,53){\vector(+1, 0){66}}
   \put( 51,58){$\Ass^\omega_{GM}$}
   \put( 33, 3){\vector(+1, 0){66}}
   \put( 51, 8){$\Ass^\omega_{GM}$}
   \put( 13,46){\vector( 0,-1){35}}
   \put(  2,25){$\S$}
   \put(114,46){\vector( 0,-1){35}}
   \put(118,25){$\S$}
  \end{picture}}}
 \end{equation}
 Classically the vector bundles of the form $\S\,TM$ on a manifold $M$
 with a smooth functor $\S$ are called pseudotensor bundles, their sections
 pseudotensors, and they comprise exactly the natural vector bundles of
 order one. Some modern authors however seem to confuse the classical concept
 of tensors with the property of having a value defined at every point.

 \begin{Lemma}[Properties of Association Functors]
 \hfill\label{assprop}\break
  Consider a principal $G$--bundle $GM$ over a manifold $M$ endowed with a
  principal connection $\omega$ and the corresponding association functor
  from the category $\MF_G$ of manifolds endowed with smooth $G$--actions
  to the category $\FB^\nabla_M$ of fiber bundles over $M$ endowed with
  non--linear connections:
  \begin{enumerate}
   \item The association functor $\Ass^\omega_{GM}$ preserves Cartesian
         products:
         $$
          GM\;\times_G\,(\,\F\,\times\,\hat\F\,)
          \;\;=\;\;
          (\;GM\,\times_G\F\;)\,\times_M\,(\;GM\,\times_G\hat\F\;)\ .
         $$
   \item On the full subcategory $\MF\,\subset\,\MF_G$ of manifolds
         with trivial $G$--action the association functor $\Ass^\omega_{GM}$
         agrees with the product functor:
         $$
          \left.\Ass^\omega_{GM}\right|_\MF:\;\;\MF\;\longrightarrow\;
          \FB^\nabla_M,\qquad\F\;\longmapsto\;M\,\times\,\F\ .
         $$
   \item Restricted to the subcategory $\Rep_G\,\subset\,\MF_G$ of finite
         dimensional smooth representations of the Lie group $G$ under
         $G$--equivariant linear maps the association functor
         $\Ass^\omega_{GM}$ takes values in the subcategory $\VB^\nabla_M$
         of vector bundles over $M$ endowed with linear connections:
         $$
          \left.\Ass^\omega_{GM}\right|_{\Rep_G}:\;\;\Rep_G\;\longrightarrow\;
          \VB^\nabla_M,\qquad V\;\longmapsto\;GM\,\times_GV\ .
         $$
  \end{enumerate}
 \end{Lemma}

 \begin{Proof}
 Of course all three statements of this lemma are easily proved directly by
 unwrapping all the definitions made above; the second statement for example
 is an elaborate description of the trivial fiber bundle isomorphism
 $$
  GM\,\times_G\F
  \;\stackrel=\longrightarrow\;GM/_{\displaystyle G}\,\times\,\F
  \;\stackrel\cong\longrightarrow\;M\,\times\,\F\ ,
 $$
 whenever $G$ acts trivially on $\F$ and thus effectively only on the first
 factor of $GM\times\F$ in the construction of the quotient $GM\times_G\F$.
 This fiber bundle isomorphism is evidently natural, it is compatible with
 all the fiber bundle homomorphisms induced by smooth maps $\varphi:\,\F
 \longrightarrow\hat\F$ between manifolds $\F$ and $\hat\F$ with trivial
 $G$--action.

 Nevertheless we think the lemma is quite interesting, because the third is
 actually a consequence of the first two statements. Combining the existence
 of additive inverses and the unity axiom $\forall v:\;1\cdot v\,=\,v$ into
 the axiom $\forall v:\;v+(-1)\cdot v\,=\,0$ we see that only three structure
 maps are needed to formulate all axioms for a vector space object $V$ in a
 category $\Cat$ in terms of commutative diagrams provided we have specified
 a field object $\mathbb{K}$:
 $$
  \cdot:\;\;\mathbb{K}\,\times\,V\;\longrightarrow\;V
  \qquad\quad
  +:\;\;V\,\times\,V\;\longrightarrow\;V
  \qquad\quad
  0:\;\;\{*\}\;\longrightarrow\;V\ .
 $$
 In the category $\MF_G$ for example we may take the manifold
 $\R$ with the trivial $G$--action as the field object $\mathbb{K}\,=\,
 \R^\triv$, the corresponding vector space objects are smooth representations
 of the Lie group $G$ over $\R$.

 On the other hand the first and second statement of the lemma assert that
 the association functor $\Ass^\omega_{GM}$ preserves Cartesian products
 and agrees with the product functor $M\,\times$ on the full subcategory
 $\MF\,\subset\,\MF_G$. In consequence $\Ass^\omega_{GM}$ sends terminal
 objects in $\MF_G$ to terminal objects in $\FB^\nabla_M$ and a representation
 $V$ to a fiber bundle $VM\,:=\,GM\times_GV$ with three parallel structure
 maps, the zero section $0:\,M\longrightarrow VM$ and:
 $$
  \cdot:\;\;\R\,\times\,VM\;\longrightarrow\;VM
  \qquad\qquad
  +:\;\;VM\,\times_MVM\;\longrightarrow\;VM\ .
 $$
 According to Lemma \ref{clcx} the parallelity of the scalar multiplication
 map alone suffices to force the non--linear connection $\P^\nabla$ on
 $VM\,\in\,\Obj\;\FB^\nabla_M$ to be a linear connection in the sense of
 Definition \ref{linx}. 
 \end{Proof}

 \pfill
 Historically the concept of principal bundles and principal connections
 arose from Cartan's beautiful idea of moving frames, which asserts that
 every vector bundle $VM$ with connection $\nabla$ lies in the image of
 the association functor $\Ass^\omega_{GM}$ for some principal bundle with
 connection. A suitable choice for the principal bundle $GM$ is the frame
 bundle with model vector space $V$
 $$
  \GL(\,M,\,VM\,)
  \;\;:=\;\;
  \{\;\;(\,p,\,F\,)\;\;|\;\;p\,\in\,M\textrm{\ and\ }F:\,
  V\stackrel\cong\longrightarrow V_pM\;\;\}\ ,
 $$
 which is a principal $\GL\,V$--bundle over $M$ with right
 multiplication given by precomposition $(p,F)\,\gamma\,=\,
 (p,F\circ\gamma)$. The tautological diffeomorphism
 $$
  \GL(\,M,\,VM\,)\;\times_{\GL\,V}V\;\stackrel\cong\longrightarrow\;VM,
  \qquad[\;(\,p,\,F\,),\;v\;]\;\longmapsto\;Fv\ ,
 $$
 is a parallel isomorphism for the principal connection on $\GL(M,VM)$
 $$
  \omega\Big(\;\left.\frac d{dt}\right|_0(\,p_t,\,F_t\,)\;\Big)
  \;\;:=\;\;
  \left.\frac d{dt}\right|_0
  F_0^{-1}\circ(\,\PT^\nabla_t\,)^{-1}\circ F_t
  \;\;\in\;\;
  \End\;V
 $$
 constructed from the parallel transport $\PT^\nabla_t:\,
 V_{p_0}M\longrightarrow V_{p_t}M$ with respect to $\nabla$ along the
 curve $t\longmapsto p_t$; the principal connection axiom of Lemma \ref{cpb}
 is particularly easy to verify using this definition for $\omega$.

 In consequence of this moving frames argument it does not make too much
 sense to ask, whether or not a vector bundle with connection is in the
 image of {\em some} association functor. The appropriate answer to this
 question for an association functor fixed in advance is definitely more
 interesting and was given in the master thesis of one of the authors.
 A closely related concept is the concept of geometric vector bundles
 defined in \cite{sw}:
 
 \begin{Proposition}[Images of Association Functors]
 \hfill\label{chimg}\break
  Let $G$ be a simply connected Lie group and let $GM$ be a principal
  $G$--bundle over a simply connected manifold $M$ endowed with a principal
  connection $\omega$. A vector bundle $VM$ with a linear connection
  $\P^\nabla$ over $M$ is isomorphic in the vector bundle category
  $\VB^\nabla_M$ to a vector bundle in the image of the association
  functor $\Ass^\omega_{GM}$, if and only if there exists a parallel
  bilinear map
  $$
   \star_\inf\,:\;\;\aut\;GM\;\times_M\,VM\;\longrightarrow\;VM,
   \qquad(\,X,\,v\,)\;\longmapsto\;X\,\star\,v\ ,
  $$
  which is a representation of the Lie algebra $\aut_pGM$ at every
  $p\,\in\,M$
  $$
   (\,\star_\inf\,)_p:\;\;\aut_pGM\;\times\;V_pM\;\longrightarrow\;V_pM
  $$
  with the additional property that the curvature of the given connection
  $\nabla$ agrees with the pointwise action of the curvature
  $R^\omega\,\in\,\Omega^2(\,M,\,\aut\,GM\,)$:
  $$
   R^\nabla_{X,\,Y}v
   \;\;=\;\;
   R^\omega_{X,\,Y}\,\star_\inf\,v\ .
  $$
 \end{Proposition}

 \begin{Proof}
 Consider to begin with the vector bundle $VM\,:=\,GM\times_GV$ associated
 to a representation $V$ of the Lie group $G$. According to our discussion of
 the infinitesimal action following Definition \ref{infact} the composition
 $$
  \star_\inf:\;\;\g\;\times\;V
  \;\stackrel{\star_\inf}\longrightarrow\;TV
  \;\stackrel\cong\longrightarrow\;V\,\times\,V
  \;\stackrel{\pr_R}\longrightarrow\;V
 $$
 is $G$--equivariant and thus gives rise to a parallel $\R$--bilinear map,
 which is a representation $(\,\star_\inf\,)_p$ of the Lie algebra $\aut_pGM$
 on $V_pM$ in every point:
 $$
  \star_\inf:\;\;\aut\;GM\,\times_M\,VM\;\longrightarrow\;VM\ .
 $$
 Conversely assume that $\star_\inf:\,\aut\,GM\times_MVM\longrightarrow VM$
 is a parallel representation of the Lie algebra bundle $\aut\,GM$ on a vector
 bundle $VM$ with a linear connection $\P^\nabla$. According to equation
 (\ref{giso}) the fiber Lie group $\Aut_pGM$ is isomorphic to $G$ in every
 point $p\,\in\,M$ and so simply connected, in consequence the infinitesimal
 action $(\,\star_\inf\,)_p$ of its Lie algebra $\aut_pGM$ integrates to a
 representation of the Lie group $\Aut_pGM$ on the vector space $V_pM$.
 Though slightly technical it is straightforward to prove that the integrated
 representation depends smoothly on the point $p\,\in\,M$
 \begin{equation}\label{intrep}
  \star:\;\;\Aut\,GM\;\times_M\,VM\;\longrightarrow\;VM\ ,
 \end{equation}
 the details of this argument are left to the reader. In addition to the
 vector bundle $VM$ with its connection $\P^\nabla$ we consider the vector
 bundle $GM\times_GV$ associated to some representation $V$ of $G$ endowed
 with the linear connection $\P^\omega$ induced by the principal connection
 $\omega$ in Proposition \ref{asscon}. The two connections determine a linear
 connection $\P^{(\omega,\nabla)}$ on the vector bundle $\Hom(\,GM\times_GV,
 \,VM\,)$ characterized by the fact that its parallel transport
 $$
  \PT^{(\omega,\nabla)}_t:\;\;\Hom(\,G_{p_0}M\times_GV,\,V_{p_0}M\,)
  \;\longrightarrow\;\Hom(\,G_{p_t}M\times_GV,\,V_{p_t}M\,)
 $$
 along an arbitrary curve $t\longmapsto p_t$ makes the following diagram
 commute
 \begin{equation}\label{ptx}
  \vcenter{\hbox{\begin{picture}(195,65)
   \put(  0,52){$G_{p_0}M\times_GV$}
   \put(  0, 0){$G_{p_t}M\times_GV$}
   \put(155,52){$V_{p_0}M$}
   \put(155, 0){$V_{p_t}M$}
   \put( 63,55){\vector(+1, 0){89}}
   \put(104,58){$F$}
   \put( 63, 3){\vector(+1, 0){89}}
   \put( 84, 7){$\PT^{(\omega,\nabla)}_tF$}
   \put( 38,48){\vector( 0,-1){38}}
   \put(  9,26){$\PT^\omega_t$}
   \put(170,48){\vector( 0,-1){38}}
   \put(174,26){$\PT^\nabla_t$}
  \end{picture}}}
 \end{equation}
 for all linear maps $F:\,G_{p_0}M\times_GV\longrightarrow V_{p_0}M$, where
 $\PT^\omega_t$ and $\PT^\nabla_t$ are the parallel transports along the same
 curve with respect to $\P^\omega$ and $\P^\nabla$.

 \pfill
 The principal idea of the proof is now to construct a parallel and actually
 flat vector subbundle of the vector bundle $\Hom(\,GM\times_GV,\,VM\,)$
 over $M$. For this purpose we consider the family of vector subspaces
 of the fibers
 \begin{eqnarray*}
  \lefteqn{\Big[\,\Hom_{\Aut\,GM}(\,GM\times_GV,\,VM\,)\,\Big]_p}
  &&
  \\
  &:=&
  \{\;\;F:\,G_pM\times_GV\longrightarrow V_pM\;\;|\;\;
  \textrm{linear and\ }\Aut_pGM\textrm{\ equivariant}\;\;\}
 \end{eqnarray*}
 of the vector bundle $\Hom(\,GM\times_GV,\,VM\,)$ in each point $p\,\in\,M$.
 In order to show that this family of subspaces is the family of fibers
 of a vector subbundle of $\Hom(\,GM\times_GV,\,VM\,)$ we observe that
 the parallel transport
 $$
  \PT^\omega_t:\;\;G_{p_0}M\,\times_GV\;\stackrel\cong\longrightarrow
  G_{p_t}M\,\times_GV
  \qquad\quad
  \PT^\nabla_t:\;\;V_{p_0}M\;\stackrel\cong\longrightarrow\;V_{p_t}M
 $$
 in both vector bundles $GM\times_GV$ and $VM$ along a curve $t\longmapsto p_t$
 is equivariant over the parallel transport with respect to the Lie group
 connection $\P^\omega$ on the automorphism bundle $\Aut\,GM$ induced by
 $\omega$. More precisely we find 
 $$
  \PT^\nabla_t\big(\;(\,p_0,\,\psi\,)\;\star\;v\;\big)
  \;\;=\;\;
  \PT^\omega_t(\,p_0,\,\psi\,)\;\star\;\PT^\nabla_tv
 $$
 for the vector bundle $VM$, because $\star_\inf:\,\aut\,GM\,\times_MVM
 \longrightarrow VM$ is parallel by assumption. In consequence the
 parallel transport $\PT^{(\omega,\nabla)}$ with respect to the linear
 connection $\P^{(\omega,\nabla)}$ specified in diagram (\ref{ptx})
 induces for all $t\,\in\,\R$ vector space isomorphisms $F\longmapsto
 \PT^\nabla_t\circ F\circ(\,\PT^\omega_t\,)^{-1}$ between:
 $$
  \Big[\,\Hom_{\Aut\,GM}(\,GM\times_GV,VM\,)\,\Big]_{p_0}
  \!\stackrel\cong\longrightarrow
  \Big[\,\Hom_{\Aut\,GM}(\,GM\times_GV,VM\,)\,\Big]_{p_t}\ .
 $$
 By assumption the underlying manifold $M$ is (simply) connected, and hence
 all vector subspaces $[\,\Hom_{\Aut\,GM}(\,GM\times_GV,\,VM\,)\,]_p$ have
 the same dimension. With parallel transport depending smoothly on the curve
 we conclude that $\Hom_{\Aut\,GM}(\,GM\times_GV,\,VM\,)$ is a genuine
 vector subbundle of $\Hom(\,GM\times_GV,\,VM\,)$, moreover it is a parallel
 subbundle as it is invariant under parallel transport along arbitrary
 curves.

 On the other hand the curvature of the linear connection
 $\P^{(\omega,\nabla)}$ on the vector bundle $\Hom(\,GM\times_GV,\,VM\,)$
 is determined by the universality of principal curvature discussed in
 Proposition \ref{asscon}, namely it holds true that
 $$
  R^{(\,\omega,\,\nabla\,)}_{X,\,Y}F
  \;\;=\;\;
  R^\nabla_{X,\,Y}\,\circ\,F\;-\;F\,\circ\,(\,R^\omega_{X,\,Y}\,\star_\inf\,)
 $$
 for all tangent vectors $X,\,Y\,\in\,T_pM$ and $F\,\in\,\Hom_p(\,GM\times_GV,
 \,VM\,)$. Due to equivariance the curvature of the connection $\P^{(\omega,
 \nabla)}$ restricted to the parallel vector subbbundle $\Hom_{\Aut\,GM}
 (\,GM\times_GV,\,VM\,)$ vanishes identically, put differently $\Hom_{\Aut\,GM}
 (\,GM\times_GV,\,VM\,)$ is a flat vector bundle over $M$ under the restriction
 of the connection $\P^{(\omega,\nabla)}$.
 
 In the argument presented so far the actual choice of the representation
 $V$ did not play any role. In order to make a diligent choice we fix a
 frame $g\,\in\,G_pM$ over a point $p\,\in\,M$ and consider the Lie group
 isomorphism
 $$
  \Phi:\;\;G\;\stackrel\cong\longrightarrow\;\Aut_pGM,\qquad\gamma\;
  \longmapsto\;\big(\,p,\,\hat g\,\longmapsto\,g\gamma(g^{-1}\hat g)\,\big)\ ,
 $$
 which is essentially the Lie group bundle isomorphism (\ref{giso}) restricted
 to the fiber of $p$. This Lie group isomorphism allows us to pull back
 the integrated representation (\ref{intrep}) of $\Aut_pGM$ on the vector
 space $V\,:=\,V_pM$ to a smooth representation $\star:\,G\times V
 \longrightarrow V$ enjoying the critical property that
 $$
  \overline\Phi:\;\;G_pM\,\times_GV\;\stackrel\cong\longrightarrow\;V
  \;\stackrel=\longrightarrow\;V_pM,\qquad [\,\hat g,\,v\,]\;\longmapsto
  \;\Phi(\,g^{-1}\hat g\,)\,\star\,v
 $$
 is an equivariant vector space isomorphism under $\Aut_pGM$ in the sense:
 \begin{eqnarray*}
  \overline\Phi\big(\;\Phi(\,\gamma\,)\,[\,\hat g,\,v\,]\;\big)
  &=&
  \overline\Phi\big(\;[\,g\gamma(g^{-1}\hat g),\,v\,]\;\big)
  \\
  &=&
  \Phi\big(\,g^{-1}g\gamma(g^{-1}\hat g)\,\big)\;\star\;v
  \;\;=\;\;
  \Phi(\,\gamma\,)\;\overline\Phi\big(\;[\,\hat g,\,v\,]\;\big)\ .
 \end{eqnarray*}
 In consequence the fiber of the vector bundle $\Hom_{\Aut\,GM}(\,GM\times_GV,
 \,VM\,)$ over the chosen point $p\,\in\,M$ contains the vector space
 isomorphism $\overline\Phi$, which translates under parallel transport along
 arbitrary curves with respect to the flat connection $\P^{(\omega,\nabla)}$
 into a parallel, globally defined section $\overline\Phi$ on the simply
 connected manifold $M$. Evaluation of this parallel section in the points
 of $M$ converts it into a parallel isomorphism of vector bundles:
 $$
  \overline\Phi:\;\;GM\,\times_GV\;\stackrel\cong\longrightarrow\;VM,\qquad
  [\,\hat g,\,v\,]\;\longmapsto\;\overline\Phi_{\pi(\hat g)}[\,\hat g,\,v\,]\ .
 $$
 \vskip-24pt
 \end{Proof}
\section{The Category of Gauge Theory Sectors}
\label{agauge}
 Every association functor is in a sense a reproducing functor, there exists
 in its source category an object, whose image in its target category is
 isomorphic to the principal bundle defining the association functor in the
 first place. Based on this simple observation we characterize the association
 functors among all functors from $\MF_G$ to $\FB^\nabla_M$ in this section,
 moreover we establish an equivalence of categories between the category of
 principal bundles and a suitably defined category of functors called gauge
 theory sectors.

 \pfill
 Consider the smooth action of a given Lie group $G$ on its underlying
 manifold by left multiplication $\star:\,G\times G\longrightarrow G,\,
 (\gamma,g)\longmapsto\gamma g$, which defines an object $G^\lft\,\in\,
 \Obj\,\MF_G$ in the category of $G$--manifolds. The image of $G^\lft$
 under the functor $\Ass^\omega_{GM}$ is isomorphic as a fiber bundle
 to $GM$
 \begin{equation}\label{reprod}
  \Ass^\omega_{GM}(\,G^\lft\,)\;\stackrel\cong\longrightarrow\;GM,
  \qquad[\,g,\,\gamma\,]\;\longmapsto\;g\,\gamma\ ,
 \end{equation}
 and the inverse isomorphism $g\longmapsto[g,e]$ is easily verified
 to be parallel with
 $$
  \P^\nabla\Big(\;\left.\frac d{dt}\right|_0[\,g_t,\,e\,]\;\Big)
  \;\;=\;\;
  \Big[\;g_0,\;\left.\frac d{dt}\right|_0e\;+\;
  \omega\Big(\,\left.\frac d{dt}\right|_0g_t\,\Big)\,\star_\inf\,e\,\Big]
  \;\;=\;\;
  0
 $$
 whenever $\left.\frac d{dt}\right|_0g_t$ is horizontal in the sense
 $\omega(\left.\frac d{dt}\right|_0g_t)\,=\,0$. This reproducing property
 of $\Ass^\omega_{GM}$ lies at the heart of the proof of the following theorem:
 
 \begin{Theorem}[Characterization of Association Functors]
 \hfill\label{natra}\break
  Consider a covariant functor $\FF:\,\MF_G\longrightarrow\FB^\nabla_M$
  from the category of $G$--manifolds to the category of fiber bundles
  with connection over $M$. If the functor $\FF$ preserves Cartesian
  products and agrees with the product functor
  $$
   M\,\times:\;\;\MF\;\longrightarrow\;\FB^\nabla_M,\qquad
   \F\;\longmapsto\;M\,\times\,\F\ ,
  $$
  on the full subcategory $\MF\,\subset\,\MF_G$ of manifolds with trivial
  $G$--action, then $\FF$ is naturally isomorphic to the association functor
  corresponding to some principal $G$--bundle $GM$ endowed with a principal
  connection $\omega$.
 \end{Theorem}
 
 \begin{Proof}
 Consider a functor $\FF:\,\MF_G\longrightarrow\FB^\nabla_M$ from the category
 of $G$--manifolds to the category of fiber bundles over $M$ endowed with
 non--linear connections, which preserves Cartesian products and agrees with
 the product functor $M\times:\,\MF\longrightarrow\FB^\nabla_M$ on the full
 subcategory of trivial $G$--manifolds. At least three different objects in
 the domain category $\MF_G$ of the functor $\FF$ have underlying manifold
 equal to the Lie group $G$:
 $$
  G^\lft
  \qquad\qquad
  G^\ad
  \qquad\qquad
  G^\triv\ .
 $$
 The difference between these three objects in $\MF_G$ resides in their
 actions, which is by left multiplication $\gamma\star g\,:=\,\gamma g$
 and conjugation $\gamma\star g\,:=\,\gamma g\gamma^{-1}$ respectively for
 $G^\lft$ and $G^\ad$, whereas $G$ acts trivially on $G^\triv$. Every terminal
 object in the category $\MF_G$ is a zero--dimensional manifold point $\{*\}$
 with necessarily trivial $G$--action, hence $\FF$ maps it to the terminal
 object $M\times\{*\}$ in the category $\FB^\nabla_M$. In other words the
 functor $\FF$ maps terminal objects to terminal objects and preserves
 Cartesian products and in consequence turns group like and principal
 objects in the category $\MF_G$ into group like and principal objects
 in the category $\FB^\nabla_M$.

 With $G$ acting by automorphisms on both $G^\ad$ and $G^\triv$ both objects
 are group like objects in the category $\MF_G$ under the multiplication
 and inverse inherited from $G$. The significance of the group like object
 $\FF(\,G^\ad\,)$ in the category $\FB^\nabla_M$ may be somewhat obscure at
 this point, the group like object $\FF(\,G^\triv\,)\,=\,M\times G$ however is
 just the trivial $G$--bundle over $M$ endowed with the trivial connection.
 Moreover the original Lie group multiplication defines $G$--equivariant
 structure maps in analogy to definition (\ref{stm})
 $$
  \rho:\;\;G^\lft\;\times\;G^\triv\;\longrightarrow\;G^\lft
  \qquad\qquad
  \backslash:\;\;G^\lft\;\times\;G^\lft\;\longrightarrow\;G^\triv
 $$
 by means of $\rho(g,\hat g)\,:=\,g\hat g$ and $\backslash(g,\hat g)\,:=\,g^{-1}
 \hat g$, which naturally enough turn $G^\lft$ into a $G^\triv$--principal
 object in the category $\MF_G$. According to Lemma \ref{glo} the image of
 $G^\lft$ is a principal $G$--bundle $GM\,:=\,\FF(\,G^\lft\,)$ over the
 manifold $M$ endowed with a principal connection $\omega$. In passing we
 observe that the group like object $G^\ad$ acts $G$--equivariantly on
 $G^\lft$ via
 $$
  \star:\;\;G^\ad\;\times\;G^\lft\;\longrightarrow\;G^\lft,
  \qquad (\,\gamma,\,g\,)\;\longmapsto\;\gamma\,g\ ,
 $$
 and this action identifies the group like object $\FF(\,G^\ad\,)$ in the
 category $\FB^\nabla_M$ with the gauge group bundle $\Aut\,GM$ of $GM$ by
 means of the action:
 $$
  \FF(\,\star\,):\;\;\FF(\,G^\ad\,)\;\times_M\,GM\;\longrightarrow\;GM\ .
 $$
 It remains to show that the original functor $\FF$ is naturally isomorphic
 to the association functor $\Ass^\omega_{GM}$. For this purpose we consider
 a general object $\F\,\in\,\Obj\;\MF_G$; replacing its $G$--action by the
 trivial $G$--action on the same underlying manifold we project it to an
 object $\F^\triv\,\in\,\Obj\;\MF$ in the subcategory of manifolds with
 trivial $G$--action. The $G$--equivariant map
 $$
  \Psi:\;\;G^\lft\,\times\,\F\;\stackrel\cong\longrightarrow\;
  G^\lft\,\times\,\F^\triv,\qquad(\,g,\,f\,)\;\longmapsto\;
  (\,g,\,g^{-1}\,\star\,f\,)
 $$
 is actually an isomorphism in $\MF_G$ with inverse $(\,g,\,f\,)\,\longmapsto
 \,(\,g,\,g\star f\,)$, which fits for an arbitrary element $\gamma\,\in\,G$
 into the commutative diagram
 \begin{equation}\label{tact}
  \vcenter{\hbox{\begin{picture}(170,100)(0,0)
   \put(  0, 0){$G^\lft\times\F$}
   \put(  0,88){$G^\lft\times\F$}
   \put(110, 0){$G^\lft\times\F^\triv$}
   \put(110,88){$G^\lft\times\F^\triv$}
   \put( 78,44){$\F$}
   \put( 52,91){\vector(+1, 0){54}}
   \put( 78,95){$\scriptstyle\Psi$}
   \put( 46,85){\vector(+1,-1){31}}
   \put( 43,67){$\scriptstyle\pr_R$}
   \put(122,85){\vector(-1,-1){31}}
   \put(110,67){$\scriptstyle\star$}
   \put( 52, 3){\vector(+1, 0){54}}
   \put( 78, 7){$\scriptstyle\Psi$}
   \put( 46,11){\vector(+1,+1){31}}
   \put( 43,28){$\scriptstyle\pr_R$}
   \put(122,11){\vector(-1,+1){31}}
   \put(110,28){$\scriptstyle\star$}
   \put( 29,84){\vector( 0,-1){72}}
   \put(  2,45){$\scriptstyle\rho_\gamma\times\id$}
   \put(139,84){\vector( 0,-1){72}}
   \put(143,45){$\scriptstyle\rho_\gamma\times(\gamma^{-1}\star)$}
  \end{picture}}}
 \end{equation}
 in the category $\MF_G$, where $\rho_\gamma:\,G^\lft\longrightarrow G^\lft,
 \,g\longmapsto g\gamma,$ denotes the right multiplication by $\gamma$ and
 $\star$ the original $G$--action characterizing the object $\F$ thought of as
 a $G$--equivariant (sic!) map $\star:\,G^\lft\times\F^\triv\longrightarrow\F$.
 Writing the right multiplication $\rho_\gamma$ in the category $\MF_G$ as
 a composition
 $$
  G^\lft
  \;\stackrel{\id\times\term}\longrightarrow\;
  G^\lft\,\times\,\{*\}
  \;\stackrel{\id\times\gamma}\longrightarrow\;
  G^\lft\,\times\,G^\triv
  \;\stackrel\rho\longrightarrow\;
  G^\lft
 $$
 factorizing over the element morphism $\gamma:\,\{*\}\longrightarrow
 G^\triv$ in the subcategory $\MF\,\subset\,\MF_G$ we conclude that
 $\FF(\,\rho_\gamma\,):\,GM\longrightarrow GM$ agrees with the right
 multiplication $R_\gamma:\,GM\longrightarrow GM,\,g\longmapsto g\gamma,$
 in the principal bundle $GM$ induced by $\FF(\,\rho\,):\,GM\times G
 \longrightarrow GM$, because $\FF$ preserves Cartesian products and
 agrees with the product functor $M\times$ on the trivial $G$--manifolds
 $\{*\}$ and $G^\triv$. In consequence the commutative diagram
 (\ref{tact}) translates under the functor $\FF$ into the following
 commutative diagram
 \begin{equation}\label{ffact}
  \vcenter{\hbox{\begin{picture}(200,100)(0,0)
   \put(  0, 0){$GM\times_M\F M$}
   \put(  0,88){$GM\times_M\F M$}
   \put(140, 0){$GM\times\F$}
   \put(140,88){$GM\times\F$}
   \put( 87,44){$\F M$}
   \put( 74,91){\vector(+1, 0){62}}
   \put(108,95){$\scriptstyle\Psi$}
   \put( 56,85){\vector(+1,-1){31}}
   \put( 53,67){$\scriptstyle\pr_R$}
   \put(144,85){\vector(-1,-1){31}}
   \put(134,67){$\scriptstyle\FF(\,\star\,)$}
   \put( 74, 3){\vector(+1, 0){62}}
   \put(108, 7){$\scriptstyle\Psi$}
   \put( 56,11){\vector(+1,+1){31}}
   \put( 53,28){$\scriptstyle\pr_R$}
   \put(144,11){\vector(-1,+1){31}}
   \put(134,28){$\scriptstyle\FF(\,\star\,)$}
   \put( 29,84){\vector( 0,-1){72}}
   \put(  0,45){$\scriptstyle R_\gamma\times\id$}
   \put(169,84){\vector( 0,-1){72}}
   \put(174,45){$\scriptstyle R_\gamma\times(\gamma^{-1}\star)$}
  \end{picture}}}
 \end{equation}
 in the category $\FB^\nabla_M$ with $\F M\,:=\,\FF(\,\F\,)$, because $\FF$
 preserves Cartesian products, hence preserves projections and agrees on
 manifolds with trivial $G$--action like $\F^\triv$ with the product
 functor $M\times$. The parallel homomorphism $\FF(\,\star\,):\,GM\times\F
 \longrightarrow\F M$ thus descends to the quotient
 $$
  \overline{\FF(\,\star\,)}:\;\;
  GM\,\times_G\F\;\longrightarrow\;\F M
 $$
 of $GM\times\F$ by the right $G$--action defining the associated
 fiber bundle $GM\times_G\F$, which lets $\gamma\,\in\,G$ act by
 $R_\gamma\times(\gamma^{-1}\star)$. It goes without saying that
 the projection $\pr_R:\,GM\times_M\F M\longrightarrow\F M$ factors
 through the quotient of $GM\times_M\F M$ by the right $G$--action
 on the principal bundle $GM$, the commutative diagram (\ref{ffact})
 ensures moreover that the quotient diagram
 $$
  \begin{picture}(203,62)
   \put(  0,50){$(GM/_{\displaystyle G})\times_M\F M$}
   \put(150,50){$GM\times_G\F$}
   \put(102, 0){$\F M$}
   \put( 96,54){\vector(+1, 0){52}}
   \put(120,58){$\scriptstyle\overline\Psi$}
   \put( 67,46){\vector(+1,-1){36}}
   \put( 67,25){$\scriptstyle\pr_R$}
   \put(165,46){\vector(-1,-1){36}}
   \put(155,25){$\scriptstyle\overline{\FF(\,\star\,)}$}
  \end{picture}
 $$
 still commutes. With $\pr_R:\,M\times_M\F M\stackrel\cong\longrightarrow\F M$
 and $\overline{\Psi}$ being parallel diffeomorphisms of fiber bundles with
 connections over $M$ we conclude that
 $$
  \overline{\FF(\,\star\,)}:
  \;\;GM\,\times_G\F\;\stackrel\cong\longrightarrow\;\F M
 $$
 is actually an isomorphism in the category $\FB^\nabla_M$,
 moreover the construction of this parallel fiber bundle isomorphism
 $\overline{\FF(\,\star\,)}:\,\Ass^\omega_{GM}\F\longrightarrow\FF(\,\F\,)$
 for a given object $\F\,\in\,\Obj\;\MF_G$ is natural under morphisms in
 $\MF_G$ and comprises a natural isomorphism $\overline{\FF(\,\cdot\,)}:
 \,\Ass^\omega_{GM}\longrightarrow\FF$ of functors.
 \end{Proof}

 \pfill
 In order to press the point of Theorem \ref{natra} home let us define two
 rather special categories associated to a smooth manifold $M$. Objects in
 the category $\PB^\nabla_M$ of principal bundles with connections over
 $M$ are triples $(\,G,\,GM,\,\omega\,)$ formed by a Lie group $G$ and a
 principal $G$--bundle $GM$ over $M$ endowed with a principal connection
 $\omega$. Every morphism between two such objects
 $$
  (\,\varphi_\grp,\,\varphi\,):\;\;(\;G,\;GM,\;\omega\;)
  \;\longrightarrow\;(\;\hat G,\;\hat GM,\;\hat\omega\;)
 $$
 consists of a parallel homomorphism $\varphi:\,GM\longrightarrow\hat GM$
 of fiber bundles which is $G$--equivariant over the Lie group homomorphism
 $\varphi_\grp:\,G\longrightarrow\hat G$. Objects in the category
 $\GTS^\nabla_M$ of gauge theory sectors on $M$ with connections are on the
 other hand tuples $(\,G,\,\FF\,)$ formed by a Lie group $G$ and a covariant
 functor $\FF:\,\MF_G\longrightarrow\FB^\nabla_M$ which preserves Cartesian
 products and agrees with the product functor on the full subcategory
 $\MF\,\subset\,\MF_G$ of manifolds with trivial $G$--action. In
 $\GTS^\nabla_M$ morphisms are again tuples
 $$
  (\,\varphi_\grp,\,\Phi\,):
  \;\;(\;G,\;\FF\;)\;\longrightarrow\;(\;\hat G,\;\hat\FF\;)
 $$
 consisting of a group homomorphism $\varphi_\grp:\,G\longrightarrow\hat G$
 between the two Lie groups and a natural transformation $\Phi:\,\FF\,\circ
 \,\varphi^*_\grp\longrightarrow\hat\FF$ between the two functors $\MF_{\hat G}
 \longrightarrow\FB^\nabla_M$ involved, where the action pull back functor
 $$
  \varphi^*_\grp:\;\;\MF_{\hat G}\;\longrightarrow\;\MF_G,\qquad
  (\,\hat\F,\,\star_{\hat G}\,)\;\longrightarrow\;(\,\hat\F,\,\star_G\,)
 $$
 induced by $\varphi_\grp$ lets $G$ act via $g\,\star_Gf\,:=\,\varphi_\grp(g)
 \,\star_{\hat G}f$ on a $\hat G$--manifold $\hat\F$. We want to interpret
 the construction of the association functor as a functor
 $$
  \Ass:\;\;\PB^\nabla_M\;\longrightarrow\;\GTS^\nabla_M
 $$
 with $(\,G,\,GM,\,\omega\,)\longmapsto(\,G,\,\Ass^\omega_{GM}\,)$ on objects,
 hence we still have to spe\-cify $\Ass$ on morphisms: Every morphism in the
 source category $\PB^\nabla_M$ is a parallel fiber bundle homomorphism
 $\varphi:\,GM\longrightarrow\hat GM$ equivariant over $\varphi_\grp:\,
 G\longrightarrow\hat G$, in the the target category $\GTS^\nabla_M$ such
 a morphism becomes the natural transformation $\Phi_\varphi$ defined for
 $\hat\F\,\in\,\Obj\,\MF_{\hat G}$ by:
 $$
  \Phi_\varphi(\,\hat\F\,):\;\;
  GM\,\times_G\hat\F\;\longrightarrow\;\hat GM\,\times_{\hat G}\hat\F,
  \qquad[\,g,\,\hat f\,]\;\longmapsto\;[\,\varphi(g),\,\hat f\,]\ .
 $$
 
 \begin{Corollary}[Association Functor as Equivalence of Categories]
 \hfill\label{equc}\break
  For every smooth manifold $M$ the association functor $\Ass$ provides an
  equi\-valence of categories from the category $\PB^\nabla_M$ of principal
  bundles to the category $\GTS^\nabla_M$ of gauge theory sectors over $M$
  with connections:
  $$
   \Ass:\;\;\PB^\nabla_M\;\stackrel\simeq\longrightarrow\;\GTS^\nabla_M,
   \qquad(\,G,\,GM,\,\omega\,)\;\longmapsto\;(\,G,\,\Ass^\omega_{GM}\,)\ .
  $$
  In particular two principal $G$--bundles endowed with principal connections
  on $M$ are isomorphic via a parallel, $G$--equivariant homomorphism of
  fiber bundles, if and only if their association functors are naturally
  isomorphic.
 \end{Corollary}

 \begin{Proof}
 According to Theorem \ref{natra} every gauge theory sector with connection
 $(\,G,\,\FF\,)$ is isomorphic in the category $\GTS^\nabla_M$ to an
 association functor $\Ass^\omega_{GM}$ for a suitable principal $G$--bundle
 $GM$ with a principal connection $\omega$. In order to prove Corollary
 \ref{equc} we thus need to show that the association functor $\Ass$ induces
 for two arbitrary objects in $\PB^\nabla_M$ a bijection of sets:
 \begin{eqnarray*}
  \lefteqn{\Ass:\;\;\Mor_{\PB^\nabla_M}\Big(\;(\,G,\,GM,\,\omega\,),
  \;(\,\hat G,\,\hat GM,\,\hat\omega\,)\;\Big)}\qquad\qquad\qquad
  &&
  \\[-4pt]
  &\stackrel\cong\longrightarrow&
  \Mor_{\GTS^\nabla_M}\Big(\;(\,G,\,\Ass^\omega_{GM}\,),
  \;(\,\hat G,\,\Ass^{\hat\omega}_{\hat GM}\,)\;\Big)
 \end{eqnarray*}
 Consider for this purpose a morphism $(\,\varphi_\grp,\,\Phi\,)$ in the
 category $\GTS^\nabla_M$ from the image object $(\,G,\,\Ass^\omega_{GM}\,)$
 to the image object $(\,\hat G,\,\Ass^{\hat\omega}_{\hat GM}\,)$. The natural
 transformation $\Phi$ applies to every object in $\MF_{\hat G}$, specifically
 for the object $\hat G^\lft$ describing the action of $\hat G$ on itself by
 left multiplication the natural transformation $\Phi$ provides a parallel
 homomorphism of fiber bundles
 $$
  \Phi(\,\hat G^\lft\,):\;\;
  GM\times_G\hat G\;\longrightarrow\;\hat GM\times_{\hat G}\hat G\ ,
 $$
 which we may use to define $\varphi:\,GM\longrightarrow\hat GM$ as the
 composition:
 \begin{equation}\label{ccp}
  \begin{array}{lccccccl}
   \varphi:\;&GM&\longrightarrow&
   GM\,\times_G\hat G&\stackrel{\Phi(\,\hat G^\lft\,)}\longrightarrow&
   \hat GM\,\times_{\hat G}\hat G&\stackrel\cong\longrightarrow&\hat GM
   \\
   &g&\longmapsto&[\,g,\,\hat e\,]&
   &[\,\hat g,\,\hat\gamma\,]&\longmapsto&\;\hat g\,\hat\gamma\ .
  \end{array}
 \end{equation}
 The argument we used in equation (\ref{reprod}) to show that the right hand
 side isomorphism $\hat GM_{\hat G}\hat G\longrightarrow\hat GM$ is parallel
 implies that $GM\longrightarrow GM\times_G\hat G$ is parallel as well, in
 consequence $\varphi:\,GM\longrightarrow\hat GM$ is a parallel homomorphism
 of fiber bundles.

 In order to show that $\varphi$ is equivariant over the group homomorphism
 $\varphi_\grp:\,G\longrightarrow\hat G$ we use the characteristic property
 of natural transformations like $\Phi$ for the right multiplication morphism
 $\rho_{\hat\gamma}:\,\hat G^\lft\longrightarrow\hat G^\lft,\,\hat g
 \longmapsto\hat g\hat\gamma$:
 $$
  \begin{picture}(305,68)(0,0)
   \put( 74,55){$GM\times_G\hat G$}
   \put(214,55){$\hat GM\times_{\hat G}\hat G$}
   \put( 74, 0){$GM\times_G\hat G$}
   \put(214, 0){$\hat GM\times_{\hat G}\hat G$}
   \put(132,59){\vector(+1, 0){77}}
   \put(150,62){$\Phi(\,\hat G^\lft\,)$}
   \put(132, 4){\vector(+1, 0){77}}
   \put(150, 7){$\Phi(\,\hat G^\lft\,)$}
   \put(103,50){\vector( 0,-1){39}}
   \put(  0,28){$(\Ass^\omega_{GM}\circ\varphi^*_\grp)(\,\rho_{\hat\gamma}\,)$}
   \put(243,50){\vector( 0,-1){39}}
   \put(247,28){$\Ass^{\hat\omega}_{\hat GM}(\,\rho_{\hat\gamma}\,)$}
  \end{picture}
 $$
 Of course the association functors $\Ass^\omega_{GM}\,\circ\,\varphi^*_\grp$
 and $\Ass^{\hat\omega}_{\hat GM}$ are explicitly specified on morphisms in
 Definition \ref{defass} and both vertical arrows turn out to be the right
 multiplication $[\,g,\,\hat{\scriptstyle\Gamma}\,]\longmapsto[\,g,\,
 \hat{\scriptstyle\Gamma}\hat\gamma\,]$ by $\hat\gamma\,\in\,\hat G$. In
 turn we find
 $$
  \varphi(g\gamma)
  \;\;=\;\;
  \Phi(\hat G^\lft)\,[\,g\gamma,\hat e\,]
  \;\;=\;\;
  \Phi(\hat G^\lft)\,[\,g,\varphi_\grp(\gamma)\,]
  \;\;=\;\;
  \varphi(g)\,\varphi_\grp(\gamma)
 $$
 for all $g\,\in\,GM$ and $\gamma\,\in\,G$ and conclude that $\varphi$
 is equivariant over $\varphi_\grp$. Eventually we consider for an
 arbitrary object $\hat\F\,\in\,\Obj\,\MF_{\hat G}$ the orbit map
 $\mathrm{orb}_{\hat f}:\,\hat G^\lft\longrightarrow\hat\F,\,
 \hat\gamma\longmapsto\hat\gamma\star\hat f,$ associated to an element
 $\hat f\,\in\,\hat\F$ as a morphism in the category $\MF_{\hat G}$
 with associated commutative diagram:
 $$
  \begin{picture}(323,68)(0,0)
   \put( 84,55){$GM\times_G\hat G$}
   \put(224,55){$\hat GM\times_{\hat G}\hat G$}
   \put( 84, 0){$GM\times_G\hat\F$}
   \put(224, 0){$\hat GM\times_{\hat G}\hat\F$\ .}
   \put(142,59){\vector(+1, 0){77}}
   \put(160,62){$\Phi(\,\hat G^\lft\,)$}
   \put(142, 4){\vector(+1, 0){77}}
   \put(166, 7){$\Phi(\,\hat\F\,)$}
   \put(113,50){\vector( 0,-1){39}}
   \put( -1,28){$(\Ass^\omega_{GM}\circ\varphi^*_\grp)
     (\,\mathrm{orb}_{\hat f}\,)$}
   \put(253,50){\vector( 0,-1){39}}
   \put(257,28){$\Ass^{\hat\omega}_{\hat GM}(\,\mathrm{orb}_{\hat f}\,)$}
  \end{picture}
 $$
 Definition \ref{defass} provides again an explicit description of the two
 vertical arrows and the top arrow reads $[\,g,\,\hat\gamma\,]\longmapsto
 [\,\varphi(g),\,\hat\gamma\,]$, the commutativity of the diagram thus
 implies that $\Phi(\,\hat\F\,)$ is given by $[\,g,\,\hat f\,]\longmapsto
 [\,\varphi(g),\,\hat f\,]$. In other words the two natural transforms
 $\Phi$ and $\Phi_\varphi$ agree on arbitrary objects and so the functor
 $\Ass$ is full, this is surjective on morphisms. In order to show that
 $\Ass$ is injective on morphisms or faithful the reader may simply verify
 that the equivariant map $GM\longrightarrow\hat GM$ defined in equation
 (\ref{ccp}) equals $\varphi$ in case we start with the natural transformation
 $\Phi\,=\,\Phi_\varphi$.
 \end{Proof}

 \pfill
 Mutatis mutandis the arguments presented in this section work without
 taking connections into account: A functor $\FF:\,\MF_G\longrightarrow\FB_M$
 is naturally isomorphic to the association functor $\Ass_{GM}$ for some
 principal bundle $GM$, if and only if $\FF$ preserves Cartesian products
 and agrees with the product functor $M\times:\,\MF\longrightarrow\FB_M$
 on the full subcategory of trivial $G$--manifolds. Suitably defined
 categories of principal bundles and gauge theory sectors then turn the
 association functor into an equivalence of categories:
 $$
  \Ass:\;\;\PB_M\;\stackrel\simeq\longrightarrow\;\GTS_M,
  \qquad(\,G,\,GM\,)\;\longmapsto\;(\,G,\,\Ass_{GM}\,)\ .
 $$
\vskip10pt\noindent
\parbox{230pt}{%
 Gustavo Amilcar Salda\~na Moncada\\
 Instituto de Matem\'aticas (Ciudad de M\'exico)\\
 Universidad Nacional Aut\'onoma de M\'exico\\[3pt]
 04510 Ciudad de M\'exico, MEXIQUE.\\[3pt]
 \texttt{gamilcar@ciencias.unam.mx}}

\vskip10pt\noindent
\parbox{230pt}{%
 Gregor Weingart\\
 Instituto de Matem\'aticas (Cuernavaca)\\
 Universidad Nacional Aut\'onoma de M\'exico\\[3pt]
 62210 Cuernavaca, Morelos, MEXIQUE.\\[3pt]
 \texttt{gw@matcuer.unam.mx}}
\end{document}